\documentclass[12pt]{amsart}
\usepackage{amsmath,amssymb}

\def\cal{\mathcal}
\newtheorem{Thm}{Theorem}
\newtheorem{Lem}{Lemma}
\newtheorem{Prop}{Proposition}

\newtheorem{Rem}{Remark}

\newtheorem{Claim}{Claim}

\newcommand{\bb}{\mathbb}

\newcommand{\lp}{\langle\hspace{-1mm}\langle}
\newcommand{\rp}{\rangle\hspace{-1mm}\rangle}

\begin{document}

\title[Center Manifold Theorem and Stability]{Center Manifold Theorem and Stability for Integral Equations 
with Infinite Delay} 

\author{
Hideaki Matsunaga}

\address{Department of Mathematical Sciences, Osaka Prefecture University,
Sakai 599-8531, Japan}
\email{hideaki@ms.osakafu-u.ac.jp} 

\bigskip 

\author{Satoru Murakami}
\thanks{
The second author is partly supported by the Grant-in-Aid 
for Scientific Research (C), No.22540211, 
Japan Society for the Promotion of Science.
}
\address{Department of Applied Mathematics, Okayama University of Science,
Okayama 700-0005, Japan}
\email{murakami@youhei.xmath.ous.ac.jp}

\bigskip

\author{Yutaka Nagabuchi}
\address{Department of Applied Science, Okayama University of Science,
Okayama 700-0005, Japan}
\email{nagabuti@das.ous.ac.jp}
\author{Nguyen Van Minh}

\address{ Department of Mathematics and Philosophy, Columbus State University,
4225 University Avenue, Columbus GA 31907. USA}
\email{nguyen\_minh2@columbusstate.edu}

\begin{abstract}
The present paper deals with autonomous integral equations with infinite delay via dynamical system approach. Existence, local exponential attractivity, and other properties of center manifold are established by means of the 
variation-of-constants formula in the phase space that is obtained in a previous paper \cite{mur}.   
Furthermore, we prove a stability reduction principle by which the stability of an autonomous 
integral equation is implied
by that of an ordinary differential equation which we call the 
"central equation".
\end{abstract}
\thanks{Corresponding author: N.V.M.; e-mail: nguyen\_minh2@columbusstate.edu
}


\date{\today}
\maketitle


\section{Introduction}

In this paper we are concerned with the integral equation with infinite 
delay  
\begin{align}
   x(t)=\int_{-\infty}^t K(t-s)x(s)ds+f(x_{t}),  \tag{E}
\end{align}
where $K$ is a measurable $m\times m$ matrix valued function with 
complex components satisfying the condition
$$
   \int_0^\infty \|K(t)\|e^{\rho t}dt <\infty 
   \quad \text{and} \quad 
   \mathrm{ess}\,\sup\{ \|K(t)\|e^{\rho t}: t\geq 0\}<\infty, 
$$
and $f$ is a nonlinear term belonging to the space $C^1(X;{\bb C}^m)$, 
the set of all continuously (Fr\'echet) differentiable functions 
mapping $X$ into ${\bb C}^m$, with the property that  $f(0)=0$ and 
$Df(0)=0$; here, $\rho$ is a positive constant which is fixed 
throughout the paper, and $X:=L^1_{\rho}({\bb R}^-;{\bb C}^m)$, 
${\bb R}^-:=(-\infty,0]$, is a Banach space which will be introduced 
in the next section as the phase space for Eq. ($E$), and $x_t$ is 
an element in $X$ defined as $x_t(\theta)=x(t+\theta)$ for 
$\theta\in {\bb R}^-$. 
One of the purpose of the paper is to establish several results 
(the existence, the (local) exponential attractivity and so on) 
on the (local) center manifolds of the equilibrium point $0$ of Eq. ($E$). 
For several kinds of equations including ordinary differential equations, 
functional differential equations, parabolic partial differential 
equations and Volterra difference equations, the subject has been thoroughly studied. For more 
information in this direction we refer the reader to the references 
\cite{aulmin,carr, chehaltan,chilat,codlev,diek,far,farhuawu,hal, hallun, hen,hirpugshu, kri, matmur, mem, 
murnag, murmin, shub, vangil,wu} and the references therein.  Recently, 
Diekmann and Gyllenberg \cite{die-gyl} have treated Eq. ($E$)  by Adjoint Semigroup Theory,
and established several results including the principle of linearized 
stability for integral equations. Motivated by the pioneering paper \cite{die-gyl}, the authors in \cite{mur} 
have also treated integral equations with infinite delay but via a dynamical system approach, and 
established a "variation-of-constants formula" (VCF, for short) 
in the phase space. To the best of our knowledge, the existence as well as applications
of invariant manifolds, in particular center manifolds for 
Eq. ($E$) are still open questions.  It is the purpose of this paper to address these questions via the VCF established in \cite{mur}.

\medskip
Center manifolds play a crucial role in the stability analysis of 
systems around non-hyperbolic equilibria. The existence as well as the smoothness of center (center-stable) manifolds allow us to reduce the stability analysis of an original system to that of its restriction to a center (center stable) manifold. This procedure was initiated by Pliss \cite{pli}, and subsequently, becomes popular in the mathematical literature on stability and applications.  For more information on the reduction principle for ordinary differential equations in finite and infinite dimensional spaces we refer the reader to \cite{aulmin,carr,chilat,hen,malsel,minwu} and the references therein. Extensions of the reduction principle to a variety of kinds of equations, including functional differential equations could be found in \cite{aulmin,farhuawu,hallun,magrua,malsel,minwu,wu} and their references. Indeed, as stated in 
\cite[Section 10.5]{hallun}, for any functional 
differential equation (FDE) with a nonhyperbolic equilibrium point $0$ the stability of the FDE around $0$ is reduced to that of an ordinary differential equation $\dot{u}=h(u)$ 
in the neighborhood of its equilibrium point $0$. The procedure can be done based on the dynamical system restricted to a center (center stable) manifold that are assumed to exist and to be sufficiently smooth.

\medskip
We now outline the presentation of our paper.
In Section 2 we present preliminary results necessary for our 
later arguments. In Subsection 3.1 of the paper, applying the VCF in \cite{mur} 
we will prove a center manifold theorem for Eq.($E$) including the 
existence and the exponential attractivity (Theorem \ref{Theorem D1a}). 
In Subsection 3.2, introducing an 
ordinary differential equation which we call the "central equation" 
of Eq. (E), we will establish the reduction principle for integral equations (Theorem \ref{Theorem D3}) that the 
stability properties for the central equation imply that of Eq. (E) in the neighborhood of its zero solution. 
As an application of Theorem \ref{Theorem D3}  to 
stability analysis of some particular equations, we will consider 
a scalar integral equation. Indeed, by calculating the 
corresponding central equation we obtain a result 
(Proposition \ref{Proposition D8}) on the stability properties 
for the equation in the critical case. Also, in Appendix, we give a proof of the 
smoothness of the center manifold. In addition, for completeness, 
we establish the existence of other invariant manifolds (stable 
manifolds and unstable manifolds etc.) for Eq. (E). 

\vspace{5mm}


\section{Notations and Some Preparatory Results}

Let ${\bb N}$, ${\bb R}^+$, ${\bb R}^-$, ${\bb R}$ and ${\bb C}$ be 
the set of natural numbers, nonnegative real numbers, nonpositive 
real numbers, real numbers and complex numbers, respectively. 
For an $m\in {\bb N}$, we denote by ${\bb C}^m$ the space of all 
$m$-column vectors whose components are complex numbers, with the 
Euclidean norm $|\cdot|$.

Given Banach spaces $(U,\|\cdot\|_U)$ and $(V, \|\cdot\|_V)$, 
we denote by ${\cal L}(U;V)$ the space of bounded linear operators 
from $U$ to $V$ with norm
$$
   \|Q\|_{{\cal L}(U;V)}:=
     \sup\big\{\|Q(u)\|_{V}/\|u\|_{U}: \,u\in U, u\not=0 \big\}
$$ 
for $Q\in {\cal L}(U;V)$, and use the symbol ${\cal L}(U)$ in place 
of ${\cal L}(U;U)$. In particular, for an $m\times m$ matrix $M$ with 
complex components, $\|M\|$ means its operator norm 
$\|M\|_{{\cal L}({\bb C}^m)}$.

For an interval $J\subset {\bb R}$ and a Banach space $U$ we denote 
by $C(J;U)$ the space of $U$-valued continuous functions on $J$, 
and by $BC(J;U)$ its subspace of bounded continuous functions on $J$. 
We also use the notation $B_{U}(r)$ which stands for the open ball 
in $U$ at the center $0$ with radius $r>0$, that is, 
$B_{U}(r)=\{u\in U: \|u\|_{U} < r \}$.

\vskip 3mm

\subsection{Phase space and initial value problems}

Let $\rho$ be a fixed positive constant, and let $X$ be the function space  
$
   L_\rho^1({\bb R}^-;{\bb C}^m)$
   that is defined to be all equivalent classes of measurable functions
   $$
   \phi :{\bb R}^-\to {\bb C}^m: 
   \text{$\phi(\theta)e^{\rho\theta}$ is integrable on ${\bb R}^-$}. 
$$ 
Clearly, $X$ is a Banach space endowed with norm
$$
   \|\phi\|_X:=\int_{-\infty}^0 |\phi(\theta)|e^{\rho\theta}d\theta, 
   \quad \phi \in X.
$$
For any function $x:(-\infty,a)\to {\bb C}^m$ and $t<a$, we define a 
function $x_t:{\bb R}^-\to {\bb C}^m$ by $x_t(\theta):=x(t+\theta)$ 
for $\theta\in {\bb R}^-$; the function $x_t$ is called the $t$-segment 
of $x(t)$.

Consider the integral equations 
\begin{align}\label{100}
   x(t)= \int_{-\infty}^t K(t-s)x(s)ds+p(t)
\end{align}
and
\begin{align}
    x(t)= \int_{-\infty}^t K(t-s)x(s)ds+f(x_t), \tag{E}
\end{align}
where we assume, throughout the paper, that the kernel $K$ is a 
measurable $m\times m$ matrix valued function with complex components 
satisfying the conditions 
\begin{align}
   \|K\|_{1,\rho}:=\int_0^\infty \|K(t)\|e^{\rho t}dt <\infty, \quad
   \|K\|_{\infty,\rho}:= \mathrm{ess}\,\sup\{ \|K(t)\|e^{\rho t}: t\geq 0\}<\infty, \qquad
\end{align}
$p\in C({\bb R};{\bb C}^m)$ and $f:X\to {\bb C}^m$ is of class $C^1$. 
Then Eq.(\ref{100}) (resp. ($E$)) can be formulated as an abstract 
equation on the space $X$ of the form 
\begin{align}\label{127b}
  x(t)=F(t,x_t), 
\end{align}
with $F(t,\phi)=L(\phi)+p(t)$ (resp. $L(\phi)+f(\phi)$) for 
$(t,\phi)\in {\bb R}\times X$, where  
$$
   L(\phi):=\int_{-\infty}^0 K(-\theta)\phi(\theta)d\theta, \quad \phi\in X.
$$
Note that, in each case, $F(t,\phi)$ is well-defined because of   
\begin{align*}
  |L(\phi)| 
  \leq \int_{-\infty}^0 \|K(-\theta)\|e^{-\rho\theta}|\phi(\theta)|e^{\rho\theta} d\theta
  \leq \|K\|_{\infty,\rho}\|\phi\|_X.
\end{align*}
Thus, $X$ may be viewed as the phase space for Eq.'s (\ref{100}) and 
($E$); in what follows we will call $X$ the phase space.  

Now let $F:[b,\infty)\times X \to {\bb C}^m$ be any continuous function, 
and consider the equation (\ref{127b}) with the initial condition
\begin{align}\label{115}
   x_\sigma=\phi, 
   \quad \text{that is}, \quad  
   x(\sigma+\theta)=\phi(\theta) 
   \quad  \text{for} \enskip \theta \in {\bb R}^-,
\end{align}
where $(\sigma, \phi) \in [b,\infty)\times X$ is given arbitrarily. 
A function $x: (-\infty,a)\to {\bb C}^m$ is said to be a solution of 
the initial value problem (\ref{127b})-(\ref{115}) on the interval 
$(\sigma, a)$ if  $x$ satisfies the 
following conditions:
\begin{enumerate}
   \item  $x_\sigma =\phi$, that is, $x(\sigma+\theta)=\phi(\theta)$ for $\theta\in {\bb R}^-$;
   \item  $x\in L_{{\rm loc}}^1[\sigma,a)$, $x$ is locally integrable on $[\sigma,a)$;
   \item  $x(t)=F(t,x_t)$ for $t\in (\sigma, a)$. 
\end{enumerate}  
If $F(t,\phi)$ is locally Lipschitz continuous in $\phi$, by 
\cite[Proposition 1]{mur} the initial value problem (\ref{127b})-(\ref{115}) 
has a unique (local) solution, which is defined globally if, in particular, 
$F(t,\phi)$ is globally Lipschitz continuous 
in $\phi$ (\cite[Proposition 3]{mur}). 
So for any $(\sigma,\phi)\in {\bb R}\times X$ (\ref{100})-(\ref{115}) 
has a unique global solution, denoted $x(t;\sigma,\phi,p)$, which is 
called the solution of Eq.(\ref{100}) through $(\sigma,\phi)$. 
Similarly, ($E$)-(\ref{115}) has a unique (local) solution, which is 
denoted by $x(t;\sigma,\phi,f)$. Moreover we remark that if $x(t)$ is 
a solution of Eq.(\ref{127b}) on $(\sigma,a)$, then $x_t$ is an 
$X$-valued continuous function on $[\sigma,a)$ (see \cite[Lemma 1]{mur}).

Now suppose that $\phi=\psi$ in $X$, that is, $\phi(\theta)=\psi(\theta)$ 
a.e. $\theta \in {\bb R}^-$. Then by the uniqueness of solutions of 
(\ref{100})-(\ref{115}) it follows that 
$x(t;\sigma,\phi,p)=x(t;\sigma,\psi,p)$ for $t\in (\sigma,\infty)$, 
so that $x_t(\sigma,\phi,p)=x_t(\sigma,\psi,p)$ in $X$ for 
$t\in [\sigma,\infty)$. In particular, given $\sigma\in {\bb R}$, 
$x_t(\sigma,\cdot,p)$ induces a transformation on $X$ for each 
$t\in [\sigma,\infty)$; and similarly for $x_t(\sigma,\cdot,f)$ 
with $t\in [\sigma, a)$ provided that $x(t;\sigma,\phi,f)$ is the 
solution of ($E$)-(\ref{115}) on $(\sigma,a)$.    

When $J$ is an interval in ${\bb R}$, a function $\xi(t)$ is called 
a solution of Eq.\,(\ref{100}) on $J$, if $\xi_{t}\in X$ is defined 
for all $t\in J$ and if it satisfies $x(t;\sigma,\xi_{\sigma},p)=\xi(t)$ 
for all $t$ and $\sigma$ in $J$ with $t\geq \sigma$; and, similarly, 
a function $\xi(t)$ is called a solution of Eq.\,($E$) on $J$ whenever 
$\xi_{t}\in X$ for $t\in J$, and $x(t;\sigma,\xi_{\sigma},f)=\xi(t)$ 
holds for all $t$ and $\sigma$ in $J$ with $t\geq \sigma$.

\vskip 3mm

\subsection{A Variation-of-constants formula and decomposition of the phase space}
 
Now, for any $t\geq 0$ and $\phi\in X$, we define $T(t)\phi \in X$ by 
\begin{align*}
 \big[T(t)\phi\big](\theta)
 &:=   x_t(\theta;,0,\phi,0)  \\
 &\,= \left\{
      \begin{array}{ll}
        x(t+\theta;0,\phi,0),  \quad & -t<\theta\leq 0, \\
        \phi(t+\theta),    & \theta \leq -t. 
      \end{array} 
     \right.
\end{align*}
Then $T(t)$ defines a bounded linear operator on $X$. We call $T(t)$ 
the solution operator of the homogeneous integral equation 
\begin{align}\label{126c}
   x(t)= \int_{-\infty}^t K(t-s)x(s)ds.
\end{align}
$\{T(t)\}_{t\geq 0}$ is a strongly continuous semigroup of bounded 
linear operators on $X$, called a solution semigroup for Eq.(\ref{126c}).   

Given a positive integer $n$, we introduce a continuous function 
$\Gamma ^n:{\bb R}^-\to {\bb R}^+$ which is of compact support 
with ${\rm supp}\,\Gamma ^n \subset [-1/n,0]$ and satisfies 
$\int_{-\infty}^0 \Gamma ^n(\theta)d\theta=1$. Obviously, 
$\Gamma ^nx\in X$ for $x\in {\bb C}^m$ and the inequality 
$\|\Gamma ^nx\|_X\leq |x|$ holds.  

The following theorem, established in \cite{mur}, gives a representation 
formula for solutions of Eq.\,(\ref{100}) in the phase space $X$, 
which is called the variation-of-constants formula (VCF, for short) 
in the phase space and plays an essential role in the present paper.

\begin{Thm}\label{Theorem A1} \cite[Theorem 3]{mur}
The segment $x_{t}(\sigma,\phi,p)$ of the solution $x(\cdot;\sigma,\phi,p)$ 
of Eq.(\ref{100}) satisfies the following relation in $X$:
$$
   x_{t}(\sigma,\phi,p)
   =T(t-\sigma)\phi+\lim_{n\to\infty}\int_{\sigma}^{t}T(t-s)(\Gamma^{n}p(s))ds,
    \quad t\geq \sigma. 
$$
\end{Thm}

Let $\overline{X}$ be a subset of $X$ of elements $\phi\in X$ which are 
continuous on $[-\varepsilon_\phi,0]$ for some $\varepsilon_\phi>0$, 
and set
$$
   X_0:=\{\psi\in X: \text{$\psi=\phi$ a.e.~on ${\bb R}^-$ for some 
   $\phi\in \overline{X}$} \}.
$$ 
Then for any $\psi\in X_0$ we can define the value of $\psi$ at 
$\theta=0$ by 
$$
  \psi[0]:=\phi(0),
$$
where $\phi$ is an element in $\overline{X}$ satisfying $\psi=\phi$ 
a.e.~on ${\bb R}^-$. It is clear that $\psi[0]$ is well-defined, and 
$X_0$ is a normed space equipped with norm
$$
   \|\psi\|_{X_0}:=\|\psi\|_X+|\psi[0]|, \quad \psi\in X_0.
$$
By \cite[Lemma 1]{mur}, we note that the solution $x(\cdot;\sigma,\phi,p)$ 
of Eq.\,(\ref{100}) through $(\sigma,\phi)\in {\bb R}\times X$ 
satisfies $x_t(\sigma,\phi,p)\in X_0$ with 
$(x_t(\sigma,\phi,p))[0]=x(t;\sigma,\phi,p)$ for $t>\sigma$.  

The following result yields an intimate relation between solutions 
of Eq.~(\ref{100}) and $X$-valued functions satisfying an integral 
equation which arises from the variation-of-constants formula in 
the phase space.

\begin{Thm}
\label{Theorem A1a} \cite[Theorem 4]{mur}
Let $p\in C({\bb R};{\bb C}^m)$.
\begin{enumerate}
   \item If $x(t)$ is a solution of Eq.~(\ref{100}) on the entire 
         ${\bb R}$, then the $X$-valued function $\xi(t):= x_t$ 
         satisfies the relations
   \begin{enumerate}
   \item $\displaystyle{\xi(t)=T(t-\sigma)\xi(\sigma)+\lim_{n\to\infty}\int_{\sigma}^tT(t-s)(\Gamma^n p(s))ds}$, $\forall \,(t,\sigma)\in {\bb R}^2$ with $t\geq \sigma$, in $X$;
   \item $\xi\in C({\bb R};X_0).$
   \end{enumerate}
   \item Conversely, if a function $\xi:{\bb R}\to X$ satisfies the relation 
$$
   \xi(t)=T(t-\sigma)\xi(\sigma)+\lim_{n\to\infty}\int_{\sigma}^tT(t-s)(\Gamma^n p(s))ds, \quad \forall \,(t,\sigma)\in {\bb R}^2\enskip \text{with}\ t\geq \sigma,
$$
then
   \begin{enumerate}
   \item [{(c)}] $\xi\in C({\bb R};X_0);$
   \item [{(d)}] if we set
$$
   u(t)=(\xi(t))[0], \quad \forall \,t\in {\bb R},
$$
then $u\in C({\bb R};{\bb C}^m)$, $u_t=\xi(t)$ (in $X$) for any 
$t\in {\bb R}$ and $u$ is a solution of Eq.~(\ref{100}) on ${\bb R}$.
   \end{enumerate}
\end{enumerate}
\end{Thm}

Based on spectral analysis of the generator $A$ of the solution semigroup 
$\{T(t)\}_{t\geq 0}$, we also have  established the decomposition theorem 
of the phase space $X$ (\cite{mur}): Let $\sigma(A)$ and $P_\sigma(A)$ be 
the spectrum and the point spectrum of the generator $A$, respectively. 
Then the following relation holds between the spectrum of $A$ and the 
characteristic roots of Eq.\,(\ref{126c})  
$$
   \sigma (A)\cap {\bb C}_{-\rho}
   = P_\sigma(A) \cap {\bb C}_{-\rho}
   = \{\lambda \in {\bb C}_{-\rho}: \det\Delta(\lambda) =0\},
$$
where ${\bb C}_{-\rho}:=\{z\in {\bb C}: \mathrm{Re}\,z>-\rho\}$, and 
$\Delta(\lambda)$ is the characteristic operator of Eq.\,(\ref{126c}), 
that is, 
$$
   \Delta(\lambda ):=\displaystyle{ E_m-\int_0^\infty K(t)e^{-\lambda t}dt}, \ 
$$
$E_m$ being the $m\times m$-unit matrix (\cite[Proposition 4]{mur}). 
Moreover, for $\mathrm{ess}\,(A)$, the essential spectrum of $A$, 
we have the estimate 
$
 \sup\limits_{\lambda \in \mathrm{ess}\,(A)} \mathrm{Re}\,\lambda \leq -\rho
$ 
(\cite[Corollarly 2]{mur}). Now set 
$\Sigma^{u}:=\{\lambda\in \sigma(A): \mathrm {Re}\,\lambda>0\}$, 
$\Sigma^{c}:=\{\lambda\in \sigma(A): \mathrm {Re}\,\lambda=0\}$, 
and $\Sigma^s:=\sigma(A)\backslash ( \Sigma^{c}\cup \Sigma^{u})$. 
Then these observations, combined with the analyticity of 
$\det\Delta (\lambda )$ on the domain ${\bb C}_{-\rho}$, yield 
the following result. 

\begin{Thm}\label{Theorem A2}\cite[Theorem 2]{mur}
Let $\{T(t)\}_{t\geq 0}$ be the solution semigroup  of Eq.(\ref{126c}). 
Then $X$ is decomposed as a direct sum of closed subspaces $E^u$, $E^c$, 
and $E^s$
$$
       X=E^{u}\oplus E^{c}\oplus E^{s}
$$
with the following properties:
\begin{enumerate}
   \item {\rm dim}\,$(E^{u}\oplus E^{c})<\infty$,
   \item $T(t)E^{u}\subset E^{u}$, $T(t)E^{c}\subset E^{c}$, and 
         $T(t)E^{s}\subset E^{s}$ for $t\in {\bb R}^+$,
   \item $\sigma(A|_{E^u})=\Sigma^u$, $\sigma(A|_{E^c})=\Sigma^c$ 
         and $\sigma(A|_{E^s\cap {\cal D}(A)})=\Sigma^s$,
   \item $T^{u}(t):=T(t)|_{E^u}$ and $T^{c}(t):=T(t)|_{E^c}$ are 
         extendable for $t\in {\bb R}$ as groups of bounded linear 
         operators on $E^u$ and $E^c$, respectively,
   \item $T^{s}(t):=T(t)|_{E^s}$ is a strongly continuous semigroup 
         of bounded linear operators on $E^s$, and its generator is 
         identical with $A|_{E^s\cap {\cal D}(A)}$,
   \item there exist positive constants  $\alpha$, $\varepsilon$ with 
         $\alpha>\varepsilon $ and a constant $C\geq 1$ such that
\begin{alignat*}{2}
   \|T^{s}(t)\|_{{\cal L}(X)} &\leq Ce^{-\alpha t}, &\quad & t\in {\bb R}^+, \\
   \|T^{u}(t)\|_{{\cal L}(X)} &\leq Ce^{\alpha t}, & & t\in {\bb R}^-, \\
   \|T^{c}(t)\|_{{\cal L}(X)} &\leq Ce^{\varepsilon |t|}, & & t\in {\bb R}.
\end{alignat*}
\end{enumerate}
\end{Thm}

In (vi) we note that $C$ is a constant depending only on $\alpha$ and 
$\varepsilon$, and that the value of $\varepsilon>0 $ can be taken 
arbitrarily small. Also, we will use the notations $E^{cu}=E^{c}\oplus E^{u}$, 
$E^{su}=E^{s}\oplus E^{u}$ etc, and denote by $\Pi^{s}$ the projection 
from $X$ onto $E^{s}$ along $E^{cu}$, and similarly for $\Pi^{u}$, 
$\Pi^{cu}$ etc. In addition, we set
$$
    C_{1}:=\|\Pi^{s}\|_{{\cal L}(X)}+\|\Pi^{c}\|_{{\cal L}(X)} 
           +\|\Pi^{u}\|_{{\cal L}(X)}.
$$

If $f\in C^1(X;{\bb C}^m)$ satisfies $f(0)=0$ and $Df(0)=0$, 
Eq.\,(\ref{126c}) is the linearized equation of Eq.($E$) around the 
equilibrium point $0$. The equilibrium point $0$ (or the zero solution) 
of Eq.\,($E$) is said to be {\it hyperbolic} provided that 
$\Delta(\lambda )$ is invertible on the imaginary axis; 
in other words, $\Sigma^c=\emptyset$.

\vskip 5mm


\section{Center Manifold Theorem for Integral Equations}

In what follows we assume that $f\in C^1(X;{\bb C}^m)$ satisfies 
$f(0)=0$ and $Df(0)=0$. In this section we will establish the 
existence of local center manifolds of the equilibrium point $0$ 
of Eq.($E$) and study their properties. To do so, in parallel with 
Eq.($E$), we will consider a modified equation of ($E$) of the form
\begin{align}
    x(t)=\int_{-\infty}^t K(t-s)x(s)ds+f_\delta(x_t), \tag{$E_\delta$}
\end{align}
where $f_\delta$ with $\delta>0$ is a modification of the original 
nonlinear term $f$; more precisely let $\chi: {\bb R}\to [0,1]$ be a 
$C^{\infty}$-function such that $\chi(t)=1$ $(|t|\leq 2)$ and $\chi(t)=0$ 
$(|t|\geq 3)$, and define 
$$
   f_{\delta}(\phi)
   :=\chi\big(\|\Pi^{su}\phi\|_X/\delta\big)\chi\big(\|\Pi^{c}\phi\|_X/\delta\big)f(\phi), \quad  \phi\in X.
$$ 
The function $f_{\delta}: X\to {\bb C}^m$ is continuous on $X$, and 
is of class $C^1$ when restricted to the open set 
$S_{\delta}:=\big\{\phi\in X: \|\Pi^{su}\phi\|_X <\delta\big\}$ 
since we may assume that $\|\Pi^{c}\phi\|_X$ is of class $C^1$ for 
$\phi\not=0$ because of $\dim E^{c}<\infty$. Moreover, by the assumption 
$f(0)=Df(0)=0$, there exist a $\delta_{1}>0$ and a nondecreasing 
continuous function $\zeta_{\ast}:(0,\delta_{1}]\to {\bb R}^+$ such that 
$\zeta_{\ast}(+\,0)=0$, 
\begin{align}\label{351}
   \|f_{\delta}(\phi)\|_{X}\leq \delta\zeta_{\ast}(\delta)
   \quad \text{and} \quad
   \|f_{\delta}(\phi)-f_{\delta}(\psi)\|_{X}
   \leq \zeta_{\ast}(\delta)\|\phi-\psi\|_X
\end{align}
for $\phi,\psi\in X$ and $\delta\in (0,\delta_{1}]$. Indeed, we may put 
$$
   \zeta_{\ast}(\delta)
   =\big(\sup_{\|\phi\|_X\leq 3\delta}\|Df(\phi)\|_{{\cal L}(X;{\bb C}^m)}\big)
    \cdot\big(1+3\sup_{0\leq t\leq 3}|\chi'(t)|\big)
$$
(cf. \cite[Lemma 4.1]{diek}). Taking $\delta_{1}>0$ small, we may also 
assume that there exists a positive number $M_{1}(\delta_1)=: M_{1}$ 
such that
\begin{align}\label{352}
    \|D f_{\delta}(\phi)\|_{{\cal L}(X;{\bb C}^m)}\leq M_{1}, 
    \quad \phi \in S_{\delta}
\end{align}
for any $\delta\in (0,\delta_1]$. Fix a positive number $\eta$ such that 
$$
   \varepsilon <\eta < \alpha,
$$
where $\varepsilon$ and $\alpha $ are the constants in 
Theorem \ref{Theorem A2}. 

\subsection{Center Manifold and its exponential attractivity}

For the existence of center manifold for Eq.$(E_\delta)$ and its 
exponential attractivity, we have the following: 

\begin{Thm}\label{Theorem D1}
There exist a positive number $\delta$ and a $C^1$-map 
$F_{\ast,\delta}:E^{c}\to E^{su}$ with $F_{\ast,\delta}(0)=0$ 
such that the following properties hold: 
\begin{enumerate}
   \item $W^{c}_\delta:={\rm graph}\,F_{\ast,\delta}$ is tangent to 
         $E^{c}$ at zero, 
   \item $W^{c}_\delta$ is invariant for Eq.\,($E_\delta$), that is, 
         if $\xi\in W^{c}_\delta$, then $x_{t}(0,\xi,f)\in W^{c}_\delta$ 
         for $t\in {\bb R}$. 
   \item Assume moreover that $\Sigma^{u}=\emptyset$. 
Then there exists a positive constant $\beta_0$ with the property that 
if $x$ is a solution of Eq.($E_\delta$) on an interval $J=[t_0,t_1]$, 
then the inequality 
$$
    \|\Pi^s x_t-F_{\ast,\delta}(\Pi^c x_t)\|_X 
    \leq C\|\Pi^s x_{t_0}-F_{\ast,\delta}(\Pi^c x_{t_0})\|_X
         e^{-\beta_0 (t-t_0)}, \quad t\in  J
$$
holds true. In particular, if $x$ is a solution on an interval 
$[t_0, \infty)$, $x_t$ tends to $W^{c}_\delta$ exponentially as 
$t\to \infty $. 
\end{enumerate}
\end{Thm} 

As will be shown in Proposition \ref{Proposition D3} given later, 
the map $F_{\ast,\delta}:E^{c}\to E^{su}$ in the above theorem is 
globally Lipschitz continuous with the Lipschitz constant 
$L(\delta)=4C^2C_1\zeta_{\ast}(\delta)/(\alpha-\eta)$. 
Noticing that $L(\delta)\to 0$ as $\delta\to 0$, one can assume 
that the number $\delta$ satisfies $\delta\in (0,\delta_1]$ 
together with $L(\delta)\leq 1$. Let us take a small $r\in (0,\delta)$ so that 
$\|F_{\ast, \delta}(\psi)\|_X<\delta$ for any $\psi\in B_{E^c}(r)$.  Such a choice of $r$ 
is possible by the continuity of $F_{\ast, \delta}$.  Set $F_\ast:=F_{\ast,\delta}|_{B_{E^c}(r)}$ 
and consider an open neighborhood $\Omega_0$ of $0$ in $X$ defined by 
$$
   \Omega_0:=\{\phi\in X: \|\Pi^{su}\phi\|_X<\delta,\ \|\Pi^{c}\phi\|_X<r\}.
$$
Observe that $f\equiv f_{\delta}$ on $\Omega_0$. Then the following 
theorem which yields a local center manifold for Eq.\,($E$) as the 
graph of $F_{\ast}$ immediately follows from Theorem \ref{Theorem D1}.

\begin{Thm}\label{Theorem D1a}
Assume that $f\in C^{1}(X;{\bb C}^m)$ with $f(0)=Df(0)=0$. Then 
there exist positive numbers $r$, $\delta$, and a $C^1$-map 
$F_{\ast}:B_{E^{c}}(r)\to E^{su}$ with $F_{\ast}(0)=0$, together 
with an open neighborhood $\Omega_0$ of $0$ in $X$, such that the 
following properties hold: 
\begin{enumerate}
   \item $W^{c}_{{\rm loc}}(r,\delta) :={\rm graph}\,F_{\ast}$ is 
         tangent to $E^{c}$ at zero, 
   \item $W_{{\rm loc}}^{c}(r,\delta)$ is locally invariant for Eq.\,($E$), 
         that is,
   \begin{enumerate}
       \item for any $\xi\in W^{c}_{{\rm loc}}(r,\delta)$ there exists a 
             $t_{\xi}> 0$ such that 
             $x_{t}(0,\xi,f)\in W^{c}_{{\rm loc}}(r,\delta)$ for $|t|\leq t_{\xi}$,
       \item if $\xi\in W^{c}_{{\rm loc}}(r,\delta)$ and 
             $x_{t}(0,\xi,f)\in \Omega_0$  for $0\leq t\leq T$, then  
             $x_{t}(0,\xi,f)\in W^{c}_{{\rm loc}}(r,\delta)$ for 
             $0\leq t\leq T$.
   \end{enumerate}
   \item Assume moreover that $\Sigma^{u}=\emptyset$. Then there 
         exists a positive constant $\beta_0$ with the property that 
         if $x$ is a solution of Eq.($E$) on an interval $J=[t_0,t_1]$ 
         satisfying $x_t\in \Omega_0$ on $J$, then the inequality 
$$
   \|\Pi^s x_t-F_{\ast}(\Pi^c x_t)\|_X 
   \leq C\|\Pi^s x_{t_0}-F_{\ast}(\Pi^c x_{t_0})\|_X e^{-\beta_0 (t-t_0)}, 
        \quad t\in J
$$
         holds true. In particular, if the solution $x(t)$ is 
         defined on $[t_0,\infty)$ satisfying $x_t\in \Omega_0$ on 
         $[t_0,\infty)$, then $x_t$ tends to $W^{c}_{{\rm loc}}(r,\delta)$ 
         exponentially as $t\to \infty $. 
\end{enumerate}
\end{Thm} 

In what follows we will prove Theorem \ref{Theorem D1} 
by establishing several propositions. Before doing so, 
we prepare the following lemma:

\begin{Lem}\label{Lemma D1} 
Let $f_\ast\in C(X;{\bb C}^m)$, and consider the equation
\begin{align}
    x(t)=\int_{-\infty}^t K(t-s)x(s)ds+f_\ast(x_t).  \tag{$E_\ast$}
\end{align}
Moreover, let $\psi\in E^{c}$, and $\eta$ be as above. Then we have:
\begin{enumerate}
   \item If $x(t)$ is a solution of Eq.\,($E_\ast$) defined on ${\bb R}$ 
         with the properties that $\Pi^{c}x_{0}=\psi$, 
         $\sup_{t\in {\bb R}} \|x_{t}\|_X\,e^{-\eta |t|}<\infty$ and 
         $\sup_{t\in {\bb R}} |f_\ast(x_{t})|<\infty$, 
         then the $X$-valued function $u(t):=x_t$ satisfies
\begin{align*}
   u(t) &= T^c(t)\psi+ \lim_{n\to \infty} \int_{0}^{t}T^c(t-s)\Pi^{c}\Gamma^n f_\ast(u(s))ds \\
        &\quad -\lim_{n\to \infty} \int_{t}^{\infty}T^u(t-s)\Pi^{u}\Gamma^n f_\ast(u(s))ds 
        + \lim_{n\to \infty} \int_{-\infty}^{t}T^s(t-s)\Pi^{s}\Gamma^n f_\ast(u(s))ds 
\end{align*}
         for $t\in{\bb R}$, and moreover $u$ belongs to $C({\bb R}; X_0)$. 
   \item Conversely, if $y\in C({\bb R};X)$ with 
         $\sup_{t\in {\bb R}} \|y(t)\|_X \,e^{-\eta |t|}<\infty$ 
         and $\sup_{t\in {\bb R}} |f_\ast(y(t))|<\infty$ satisfies 
\begin{align*}
   y(t) &= T^c(t)\psi+\lim_{n\to \infty} \int_{0}^{t}T^c(t-s)\Pi^{c}\Gamma^n f_\ast(y(s))d\tau \\
        &\quad -\lim_{n\to \infty} \int_{t}^{\infty}T^u(t-s)\Pi^{u}\Gamma^n f_\ast(y(s))ds
        + \lim_{n\to \infty} \int_{-\infty}^{t}T^s(t-s)\Pi^{s}\Gamma^n f_\ast(y(s))ds
\end{align*}
        for $t\in {\bb R}$, then $y$ belongs to $C({\bb R};X_0)$ and the 
        function $\xi(t)$ defined by
$$
    \xi(t):= \big(y(t)\big)[0],  \quad t \in {\bb R} 
$$
        is a solution of Eq.(E$_\ast$) on ${\bb R}$ satisfying 
        $\Pi^{c}\xi_{0}=\psi$, 
        $\sup_{t\in {\bb R}} \|\xi_{t}\|_X\,e^{-\eta |t|}<\infty$ 
        and $\xi_{t}=y(t)$ for $t\in {\bb R}$.
\end{enumerate}
\end{Lem}

\noindent
{\it Proof}. (i) Let $p(t)=f_\ast(u(t))$ for $t\in {\bb R}$. Then 
$p$ belongs to $BC({\bb R};{\bb C}^m)$ and $x$ is a solution of 
Eq.\,(\ref{100}) on ${\bb R}$. So it follows from Theorem \ref{Theorem A1a} 
that $u(t)=x_t$ is a continuous $X_0$-valued function for $t\in {\bb R}$, 
that is, $u\in C({\bb R};X_0)$. We know from Theorem \ref{Theorem A1} that 
\begin{align}\label{286a}
    u(t)=T(t-\sigma)u(\sigma)+\lim_{n\to\infty}\int_{\sigma}^{t}T(t-s)\Gamma^n p(s)ds
\end{align}
holds for any $t$ and $\sigma $ with $t\geq \sigma$, which implies 
$$
 \Pi^{cu} u(t)=T^{cu}(t-\sigma)\Pi^{cu} u(\sigma)+\lim_{n\to\infty}\int_{\sigma}^{t}T^{cu}(t-s)\Pi^{cu}\Gamma^n p(s)ds, 
$$
and moreover  
\begin{align*}
   \Pi^{cu}u(\sigma)
   &= T^{cu}(\sigma-t)\Big[ \Pi^{cu}u(t)-\lim_{n\to\infty}
      \int_{\sigma}^{t}T^{cu}(t-s)\Pi^{cu}\Gamma^n p(s)ds\Big]  \\
   &= T^{cu}(\sigma-t)\Pi^{cu}u(t) - \lim_{n\to\infty}
      \int_{\sigma}^{t} T^{cu}(\sigma-s)\Pi^{cu}\Gamma^n p(s)ds, 
      \quad t\geq \sigma
\end{align*}
because of the group property of $\{T^{cu}(t)\}_{t\in {\bb R}}$. 
So it follows that 
$$
   \Pi^{cu} u(t)
   = T^{cu}(t-\sigma )\Pi^{cu} u(\sigma )+\lim_{n\to\infty}
     \int_{\sigma}^{t}T^{cu}(t-s)\Pi^{cu}\Gamma^n p(s)ds
$$
for any $t,\sigma \in {\bb R}$, and in particular 
\begin{align}\label{286a1}
   \Pi^{u} u(\sigma)
   = T^{u}(\sigma-t) \Pi^u u(t)-\lim_{n\to\infty}
     \int_{\sigma}^{t}T^{u}(\sigma-s)\Pi^{u}\Gamma^n p(s)ds, 
     \quad \forall \,t, \sigma \in {\bb R} 
\end{align}
and
\begin{align}\label{286a2}
   \Pi^{c} u(t)
   = T^{c}(t) \psi +\lim_{n\to\infty}
     \int_{0}^{t}T^{c}(t-s)\Pi^{c}\Gamma^n p(s)ds, 
     \quad t \in {\bb R}. 
\end{align}

Now, in view of  
$$
   \|T^u(\sigma-t)\Pi^{u}u(t)\|_X 
   \leq CC_{1}e^{\alpha(\sigma-t)}e^{\eta t}\Big(
        \sup_{t\in {\bb R}} \|x_{t}\|_X \,e^{-\eta t} \Big)\to 0
        \quad \text{as} \quad t\to\infty,
$$
we get from (\ref{286a1})
$$
   \Pi^{u}u(\sigma) 
   = -\lim_{t\to\infty}\lim_{n\to\infty}
     \int_{\sigma}^{t}T^u(\sigma-s)\Pi^{u}\Gamma^n p(s)ds, 
     \quad \sigma \in {\bb R}.
$$
Note that the limit 
$$
   \lim_{t\to\infty}\int_{\sigma}^{t}T^u(\sigma-s)\Pi^{u}\Gamma^n p(s)ds
   = \int_{\sigma}^{\infty}T^u(\sigma-s)\Pi^{u}\Gamma^n p(s)ds
$$
exists in $X$. Indeed, using the inequality $\|\Gamma ^{n}x\|_X \leq |x|$ 
for $x\in {\bb C}^m$, we have for $t_{2}\geq t_{1}\geq \sigma$ 
$$
   \left\| \int_{t_{1}}^{t_{2}}T^u(\sigma-s)\Pi^u\Gamma^n p(s)ds \right\|_X 
   \leq \frac{CC_{1}}{\alpha}
        \Big(\sup_{t\in {\bb R}}|p(t)| \Big)\, e^{\alpha(\sigma-t_{1})} \to 0
$$
as $t_{1}\to \infty$, which implies the existence of the limit, together 
with uniformity in $n$ of the convergence. On the other hand, since
\begin{align*}
   &\left\| \int_{\sigma}^{\infty}T^u(\sigma-s)\Pi^u\Gamma^n p(s)ds 
           -\int_{\sigma}^{\infty}T^u(\sigma-s)\Pi^u\Gamma^m p(s)ds \right\|_X    \\
   &\leq \left\| \int_t^{\infty}T^u(\sigma-s)\Pi^u\Gamma^n p(s)ds\right\|_X
         +\left\| \int_t^{\infty}T^u(\sigma-s)\Pi^u\Gamma^m p(s)ds\right\|_X \\
   &\quad +\left\| \int_{\sigma}^{t}T^u(\sigma-s)\Pi^u\Gamma^n p(s)ds 
         - \int_{\sigma}^{t}T^u(\sigma-s)\Pi^u\Gamma^m p(s)ds \,\right\|_X \\
   &\leq \frac{2CC_{1}}{\alpha}\Big(\sup_{t\in {\bb R}}|p(t)| \Big) e^{\alpha(\sigma-t)} \\
   &\quad +\left\| \Pi^u\Big( \int_{\sigma}^{t}T(\sigma-s)\Gamma^n p(s)ds 
                            - \int_{\sigma}^{t}T(\sigma-s)\Gamma^m p(s)ds \Big)\right\|_X,  
\end{align*}
for any $t>\sigma$, it follows from 
$\lim\limits_{n\to\infty}\int_\sigma ^t T(t-s)\Gamma^n p(s)ds=x_t(\sigma,0,p)$ 
(Theorem \ref{Theorem A1}) that 
\begin{align*}
  & \limsup_{n,m\to\infty}
    \left\| \int_{\sigma}^{\infty}T^u(\sigma-s)\Pi^u\Gamma^n p(s)ds 
    - \int_{\sigma}^{\infty}T^u(\sigma-s)\Pi^u\Gamma^m p(s)ds  \right\|_X \\
  &\leq \frac{2CC_{1}}{\alpha}\Big(\sup_{t\in {\bb R}}|p(t)|\Big) e^{\alpha(\sigma-t)}. 
\end{align*}
Since $t>\sigma$ is arbitrary, 
$\int_{\sigma}^{\infty}T^u(\sigma-s)\Pi^u\Gamma^n p(s)ds$ converges 
in $X$ as $n\to\infty$. These observations yield
\begin{align} \label{286b}
   \lim_{n\to\infty}\int_{\sigma}^{\infty}T^u(\sigma-s)\Pi^u\Gamma^n p(s)
   &= \lim_{n\to\infty}\lim_{t\to\infty}\int_{\sigma}^{t}T^u(\sigma-s)\Pi^u \Gamma^n p(s)ds \notag \\
   &= \lim_{t\to\infty}\lim_{n\to\infty} \int_{\sigma}^{t}T^u(\sigma-s)\Pi^u \Gamma^n p(s)ds \notag \\
   &= -\Pi^u u(\sigma), \quad \sigma \in {\bb R}. 
\end{align}
Similarly since (\ref{286a}) also implies 
$$
   \Pi^s u(t)
   =T^s(t-\sigma)\Pi^s u(\sigma)+\lim_{n\to\infty}
    \int_{\sigma}^{t}T^s(t-s)\Pi^s \Gamma^n p(s)ds, \quad t\geq \sigma, 
$$
by the same reasoning as above, one can obtain 
\begin{align}\label{286a3}
   \Pi^s u(t)=\lim_{n\to\infty}\int_{-\infty}^{t}T^s(t-s)\Pi^s \Gamma^n p(s)ds, \quad t\in {\bb R}.
\end{align}
Thus, (\ref{286a2}), (\ref{286b}) and (\ref{286a3}) yield 
\begin{align*}
   u(t)
   &= \Pi^{c}u(t)+\Pi^u u(t) +\Pi^{s}u(t) \\
   &= T^c(t)\psi+ \lim_{n\to\infty}\int_{0}^{t}T^c(t-s)\Pi^{c} \Gamma^n p(s)ds \\
   &\quad -\lim_{n\to\infty}\int_{t}^{\infty}T^u(t-s)\Pi^{u} \Gamma^n p(s)ds
          +\lim_{n\to\infty}\int_{-\infty}^{t}T^s(t-s)\Pi^s \Gamma^n p(s)ds, 
\end{align*}
$t\in {\bb R}$, as required.  

\noindent
(ii) \ Set $g(t)=f_\ast(y(t))$ for $t\in {\bb R}$. Then $g$ belongs 
to $BC({\bb R};{\bb C}^m)$, and by the same argument as the one in (i), 
the limits $\lim\limits_{n\to\infty}\int_{t}^{\infty}T^u(t-s)\Pi^u \Gamma^n g(s)ds$ and 
$\lim\limits_{n\to\infty}\int_{-\infty}^{t}T^s(t-s)\Pi^s \Gamma^n g(s)ds$ 
exist for each $t\in {\bb R}$. For any $(t,\sigma)\in {\bb R}^2$ with 
$t\geq \sigma$, we get the relation
$$
   y(t)=T(t-\sigma)y(\sigma)
        +\lim_{n\to\infty}\int_{\sigma}^tT(t-s)\Gamma^ng(s)ds
$$
in $X$, because 
{\allowdisplaybreaks
\begin{align*}
   &T(t-\sigma )y(\sigma)+ \lim_{n\to\infty}\int_{\sigma }^{t}T(t-s)\Gamma^n g(s)ds \\
   &= T(t-\sigma )\Big\{ T^c(\sigma )\psi +\lim_{n\to\infty}\int_{0}^{\sigma }T^c(\sigma -s)\Pi^c \Gamma^n g(s)ds \\
   &\quad - \lim_{n\to\infty}\int_{\sigma }^{\infty}T^u(\sigma -s)\Pi^u \Gamma^n g(s)ds 
     + \lim_{n\to\infty}\int_{-\infty}^{\sigma }T^s(\sigma -s)\Pi^s \Gamma^n g(s)ds \Big\}  \\
   &\quad + \lim_{n\to\infty}\int_{\sigma }^{t}T(t-s)\Gamma^n g(s)ds \\
   &= T^c(t)\psi + \lim_{n\to\infty}\int_{0}^{\sigma }T^c(t-s)\Pi^c \Gamma^n g(s)ds   
      - \lim_{n\to\infty}\int_{\sigma }^{\infty}T^u(t-s)\Pi^u \Gamma^n g(s)ds \\
   &\quad +\lim_{n\to\infty}\int_{-\infty}^{\sigma}T^s(t-s)\Pi^s \Gamma^n g(s)ds  
             +  \lim_{n\to\infty}\int_{\sigma }^{t}T(t-s) \Gamma^n g(s)ds \\
   &= T^c(t)\psi +\lim_{n\to\infty}\int_{0}^{t}T^c(t-s)\Pi^{c} \Gamma^n g(s)ds 
      -\lim_{n\to\infty}\int_{t}^{\infty}T^u(t-s)\Pi^u \Gamma^n g(s)ds \\
   &\quad +\lim_{n\to\infty}\int_{-\infty}^{t}T^s(t-s)\Pi^s \Gamma^n g(s)ds   \\               
   &= y(t)
\end{align*}}%
in $X$. Therefore, from Theorem \ref{Theorem A1a} it follows that 
$y\in C({\bb R};X_0)$ and the function $\xi$ defined by $\xi(t)=(y(t))[0]$, 
$t\in {\bb R}$, satisfies $\xi\in C({\bb R}; {\bb C}^m)$, $\xi_t=y(t)$ 
(in $X$) for $t\in {\bb R}$ and $\xi$ is a solution of Eq.~(\ref{100}) 
(with $p=g$) on ${\bb R}$. Observe that $\Pi^c\xi_0=\Pi^c\psi=\psi$ and 
$\sup_{t\in {\bb R}} \|\xi_{t}\|_X\,e^{-\eta |t|}= \sup_{t\in {\bb R}} \|y(t)\|_X\,e^{-\eta |t|}<\infty$. Also, since $f_\ast(\xi_{t})=f_\ast(y(t))=g(t)$ 
on ${\bb R}$, $\xi$ must be a solution of Eq.\,($E_\ast$) on ${\bb R}$. 
Thus $\xi$ is a solution of Eq.\,($E_\ast$) with the desired properties. 
The proof is completed. \qed

\vskip 5mm

Now take a $\delta_1>0$ sufficiently small so that 
\begin{align}\label{286a5}
    \zeta_{\ast}(\delta_1)CC_{1}
    \left(\frac{1}{\eta-\varepsilon }+ \frac{2}{\alpha+\eta } 
    +\frac{2}{\alpha -\eta} \right)<\frac{1}{2}
\end{align}
holds. Let $Y_{\eta}$ be the Banach space  
$$
   Y_{\eta}:=\big\{y\in C({\bb R}; X): \sup_{t\in {\bb R}} \|y(t)\|_X\,e^{-\eta|t|}<\infty\big\}
$$
with norm
$$
   \|y\|_{Y_{\eta}}
   :=\sup_{t\in {\bb R}} \|y(t)\|_X\,e^{-\eta |t|}, \quad y\in  Y_\eta,
$$
and for $(\psi, y)\in E^c \times Y_{\eta}$ set 
\begin{align}\label{287}
   {\cal F}_{\delta}(\psi,y)(t)
   &:= T^c(t)\psi + \lim_{n\to\infty}
       \int_{0}^{t}T^c(t-s)\Pi^{c}\Gamma^{n}f_{\delta}(y(s))ds \notag \\
   &\quad -\lim_{n\to\infty}\int_{t}^{\infty}
       T^u(t-s)\Pi^{u}\Gamma^{n}f_{\delta}(y(s))ds \\
   &\quad +\lim_{n\to\infty}\int_{-\infty}^{t}
       T^s(t-s)\Pi^{s}\Gamma^{n}f_{\delta}(y(s))ds \notag
\end{align}
for $t\in {\bb R}$. Notice that the right-hand side is well-defined and 
that ${\cal F}_{\delta}(\psi,y)$ is an $X$-valued function on ${\bb R}$ 
for each $(\psi, y)\in E^c \times Y_{\eta}$.   

In the following we will establish several propositions to prove 
Theorem \ref{Theorem D1}.

\begin{Prop}\label{Proposition D1} Let ${\cal F}_{\delta}(\psi, y)$ 
be as above. Then 
\begin{enumerate}
   \item ${\cal F}_{\delta}$ defines a map from $E^{c}\times Y_\eta$ 
         to $Y_\eta$ by sending $(\psi,y)\in E^{c}\times Y_\eta$ to 
         ${\cal F}_{\delta}(\psi, y)$. 
   \item Let $\delta\in (0,\delta_1]$. Then ${\cal F}_{\delta}(\psi, \cdot)$ 
         is a contraction map from $Y_{\eta}$ into itself, 
         with Lipschitz constant $1/2$, for each $\psi\in E^c$.
\end{enumerate}
\end{Prop}

\noindent
{\it Proof.} (i) We first show ${\cal F}_{\delta}(\psi,y)\in C({\bb R};X)$. 
Let 
$
   z(t):=\lim\limits_{n\to\infty}\int_{-\infty}^{t}
         T^s(t-s)\Pi^{s}\Gamma^{n}f_{\delta}(y(s))ds=\lim_{n\to\infty}\int_{-\infty}^tT^s(t-s)\Pi^s\Gamma^np_{\delta}(s)ds$ for $t\in {\bb R}, 
$
where $p_{\delta}(t):=f_{\delta}(y(t))$.  Take a $\sigma$ so that $\sigma<t$. Then 
\begin{align*}
z(t)&=\lim_{n\to\infty}\int_{-\infty}^{\sigma}T^s(t-s)\Pi^s\Gamma^np_{\delta}(s)ds+
\lim_{n\to\infty}\int_{\sigma}^{t}T^s(t-s)\Pi^s\Gamma^np_{\delta}(s)ds\\
&=T^s(t-\sigma)\bigg(\lim_{n\to\infty}\int_{-\infty}^{\sigma}T^s(\sigma-s)\Pi^s\Gamma^np_{\delta}(s)ds\bigg)+\Pi^s\bigg(
\lim_{n\to\infty}\int_{\sigma}^{t}T(t-s)\Gamma^np_{\delta}(s)ds\bigg)\\
&=\Pi^s\left\{T(t-\sigma)z(\sigma)+\lim_{n\to\infty}\int_{\sigma}^{t}T(t-s)\Gamma^np_{\delta}(s)ds\right\}.
\end{align*}
Observe that the $X$-valued function $T(t-\sigma)z(\sigma)$ is continuous in $t$. Also, the term 
$\lim_{n\to\infty}\int_{\sigma}^tT(t-s)\Gamma^np_{\delta}(s)ds$ is continuous in $t$ as an $X$-valued function, because of 
$\lim_{n\to\infty}\int_{\sigma}^tT(t-s)\Gamma^np_{\delta}(s)ds=x_t(\sigma,0,p_{\delta})$ by Theorem \ref{Theorem A1}.  
This observation leads to the continuity of $z(t)$ on ${\bb R}$.  Almost the same argument as for $z(t)$ shows the continuity of the third term 
 of the right-hand 
side of (\ref{287}).  The second term of the right-hand 
side of (\ref{287}) is identical with 
$\Pi^c x_t(0,0,p_\delta)$ for $t\in {\bb R}$ (cf. (\ref{286a2})), and 
hence it is continuous on ${\bb R}$.  Thus ${\cal F}_{\delta}(\psi,y)$ 
belongs to $C({\bb R};X)$.

By virtue of Theorem \ref{Theorem A2} combined with (\ref{351}), 
we see that 
\begin{align}\label{357}
   \big\|\big[\Pi^{c}{\cal F}_{\delta}(\psi, y)\big](t)\big\|_Xe^{-\eta |t|}
   &\leq e^{-\eta |t|} \left( Ce^{\varepsilon |t|}\|\psi\|_X 
         +\left|\int_{0}^{t}CC_{1} e^{\varepsilon |t-s|}\delta \zeta_{\ast}(\delta)ds\right| \right) \notag \\ 
   &\leq C\|\psi\|_X + \frac{ CC_{1} \delta \zeta_{\ast}(\delta )}{\varepsilon }
\end{align}
and that
\begin{align}\label{358}
   \big\|\big[\Pi^{su}{\cal F}_{\delta}(\psi, y)\big](t)\big\|_X e^{-\eta |t|}
   &\leq e^{-\eta |t|} \left(\int_{t}^{\infty}CC_{1} e^{\alpha (t-s)}\delta \zeta_{\ast}(\delta)ds 
           +\int_{-\infty}^{t}CC_{1} e^{-\alpha (t-s)}\delta \zeta_{\ast}(\delta)ds \right) \notag \\
   &\leq \left(\frac{2CC_{1} \delta \zeta_{\ast}(\delta )}{\alpha}\right)e^{-\eta |t|}
\end{align}
for $(\psi,y)\in E^c\times Y_{\eta}$ and $t\in {\bb R}$. So it follows that
$$
   \big\|\big[{\cal F}_{\delta}(\psi, y)\big](t)\big\|_X e^{-\eta |t|}
   \leq C\|\psi\|_X+\delta \zeta_{\ast}(\delta ) CC_{1}
   \left( \frac{1}{\varepsilon}+\frac{2}{\alpha} \right), 
   \quad t\in {\bb R};
$$
hence ${\cal F}_{\delta}(\psi, y)$ belongs to $Y_\eta$ with 
$
   \|{\cal F}_{\delta}(\psi, y)\|_{Y_{\eta}} 
   \leq C\|\psi\|_X+\delta\zeta_{\ast}(\delta)CC_1(1/\varepsilon+2/\alpha). 
$
Thus, ${\cal F}_{\delta}$ defines a map from $E^c\times Y_{\eta}$ 
to $Y_\eta$. 

\noindent
(ii) Let $\psi\in E^c$ and $y_{1},y_{2}\in Y_{\eta}$. Then by (\ref{351}), 
together with (\ref{286a5}),   
\begin{align*}
   \|{\cal F}_{\delta}(\psi,y_{1})-{\cal F}_{\delta}(\psi,y_{2})\|_{Y_{\eta}} 
   &\leq \sup_{t\in {\bb R}}\, e^{-\eta |t|} \left|\int_{0}^{t}
          CC_{1} \zeta_{\ast}(\delta ) e^{-\varepsilon (t-s)}\|y_{1}- y_{2}\|_{Y_{\eta}}  e^{\eta |s|}ds\right|  \\
   &\quad +\sup_{t\in {\bb R}} \,e^{-\eta |t|} \int_{t}^{\infty}
          CC_{1} \zeta_{\ast}(\delta )e^{\alpha (t-s)}\|y_{1}-y_{2}\|_{Y_{\eta}}  e^{\eta |s|}ds \\
   &\quad +\sup_{t\in {\bb R}} \,e^{-\eta |t|} \int_{-\infty}^{t}
          CC_{1} \zeta_{\ast}(\delta )e^{-\alpha (t-s)}\|y_{1}-y_{2}\|_{Y_{\eta}}  e^{\eta |s|}ds \\
  &\leq \zeta_{\ast}(\delta_1)CC_{1}\left(\frac{1}{\eta-\varepsilon }
        + \frac{2}{\alpha+\eta} +\frac{2}{\alpha -\eta} \right) 
        \|y_{1}- y_{2}\|_{Y_{\eta}} \\
  &\leq  (1/2)\|y_{1}- y_{2}\|_{Y_{\eta}},
\end{align*}
so that ${\cal F}_{\delta}(\psi, \cdot)$ is a contraction map with 
the required property.  \qed

\vskip 5mm

In view of Proposition \ref{Proposition D1} (ii) the map 
${\cal F}_{\delta}(\psi, \cdot)$ has a unique fixed point for each 
$\psi\in E^c$, say $\Lambda_{\ast,\delta}(\psi)\in Y^{\eta}$, i.e., 
we have
\begin{align}\label{363}
   \Lambda_{\ast,\delta} (\psi)(t)
   &= T^c(t)\psi + \lim_{n\to \infty} \int_{0}^{t}
      T^c(t-s)\Pi^{c}\Gamma^{n}f_{\delta}(\Lambda_{\ast,\delta}(\psi)(s))ds \notag\\
   &\quad - \lim_{n\to \infty} \int_{t}^{\infty}
      T^u(t-s)\Pi^{u}\Gamma^{n}f_{\delta}(\Lambda_{\ast,\delta}(\psi)(s))ds \\
   &\quad + \lim_{n\to \infty} \int_{-\infty}^{t}
      T^s(t-s)\Pi^{s}\Gamma^{n} f_{\delta}(\Lambda_{\ast,\delta}(\psi)(s))ds \notag
\end{align}
for $t\in {\bb R}$, whenever $0<\delta\leq \delta_1$.

\begin{Prop}\label{Proposition D2}
$\Lambda_{\ast,\delta}(\psi)$ satisfies the following:
\begin{enumerate}
   \item $\| \Lambda_{\ast,\delta}(\psi_1) -\Lambda_{\ast,\delta}(\psi_2) 
             \|_{Y_\eta}\leq 2C\|\psi_1-\psi_2\|_X$ for $\psi_1,\psi_2\in E^c$.
   \item $\Lambda_{\ast,\delta}(\psi)(t+\tau)
         =\Lambda_{\ast,\delta}\big(\Pi^c (\Lambda_{\ast,\delta}(\psi)(\tau))\big)(t)$ holds for $t,\tau\in {\bb R}$.
\end{enumerate} 
\end{Prop}

\noindent
{\it Proof.} (i) By Proposition \ref{Proposition D1} (ii) it follows that
\begin{align*}
   \|\Lambda_{\ast}(\psi_{1})- \Lambda_{\ast}(\psi_{2})\|_{Y_{\eta}}
   &= \|{\cal F}_{\delta}(\psi_{1}, \Lambda_{\ast,\delta }(\psi_{1}))
      -{\cal F}_{\delta}(\psi_{2}, \Lambda_{\ast,\delta }(\psi_{2}))\|_{Y_{\eta}} \\
   &\leq \|{\cal F}_{\delta}(\psi_{1}, \Lambda_{\ast,\delta }(\psi_{1}))
           -{\cal F}_{\delta}(\psi_{1}, \Lambda_{\ast,\delta }(\psi_{2}))\|_{Y_{\eta}} \\
   &\quad + \|{\cal F}_{\delta}(\psi_{1}, \Lambda_{\ast,\delta }(\psi_{2})) 
          -{\cal F}_{\delta}(\psi_{2}, \Lambda_{\ast,\delta }(\psi_{2}))\|_{Y_{\eta}} \\
   &\leq  (1/2) \|\Lambda_{\ast,\delta }(\psi_{1})- \Lambda_{\ast,\delta }(\psi_{2})\|_{Y_{\eta}} + \|T^{c}(\cdot) (\psi_{1}-\psi_{2})\|_{Y_\eta},
\end{align*}
so that 
\begin{align*}
   \|\Lambda_{\ast}(\psi_{1})- \Lambda_{\ast}(\psi_{2})\|_{Y_{\eta}}
   &\leq 2 \|T^{c}(\cdot) (\psi_{1}-\psi_{2})\|_{Y_\eta} \\
   &= 2 \sup_{t\in {\bb R}} \|T^c(t)(\psi_1-\psi_2)\|_X e^{-\eta|t|} \\
   &\leq 2\sup_{t\in {\bb R}} \big(Ce^{\varepsilon |t|}\|\psi_1-\psi_2\|_X\big)e^{-\eta|t|} \\
   &= 2C\|\psi_1-\psi_2\|_X.
\end{align*}

\noindent
(ii) Given $\tau\in {\bb R}$, let us set
$$
   \tilde \Lambda (t):= \Lambda_{\ast,\delta}(\psi)(t+\tau), 
   \quad t\in {\bb R}.  
$$ 
Obviously, $\tilde \Lambda (\cdot) \in Y_\eta$ and
{\allowdisplaybreaks
\begin{align*}
   \tilde \Lambda (t) 
   &= \Lambda_{\ast,\delta}(\psi)(t+\tau) \\
   &= T^c(t+\tau)\psi + \lim_{n\to \infty} \int_{0}^{t+\tau}T^c(t+\tau-s)\Pi^{c}\Gamma^{n}f_{\delta}(\Lambda_{\ast,\delta} (\psi)(s)) ds \\
   &\quad - \lim_{n\to \infty} \int_{t+\tau}^{\infty}T^u(t+\tau-s)\Pi^{u}\Gamma^{n}f_{\delta}(\Lambda_{\ast,\delta} (\psi)(s)) ds \\
   &\quad + \lim_{n\to \infty} \int_{-\infty}^{t+\tau} T^s(t+\tau-s)\Pi^{s}\Gamma^{n}f_{\delta}(\Lambda_{\ast,\delta} (\psi)(s)) ds \\
   &= T^c(t)\Big(T^c(\tau )\psi 
      + \lim_{n\to \infty} \int_{0}^{\tau}T^c(\tau-s)\Pi^{c}\Gamma^{n}f_{\delta}(\Lambda_{\ast,\delta} (\psi)(s)) ds\Big) \\
   &\quad +\lim_{n\to \infty} \int_{0}^{t}T^c(t-s)\Pi^{c}\Gamma^{n}f_{\delta}(\Lambda_{\ast,\delta} (\psi)(s+\tau)) ds \\
   &\quad -\lim_{n\to \infty} \int_{t}^{\infty}T^u(t-s)\Pi^{u}\Gamma^{n}f_{\delta}(\Lambda_{\ast,\delta} (\psi)(s+\tau)) ds \\
   &\quad +\lim_{n\to \infty} \int_{-\infty}^{t} T^s(t-s)\Pi^{s}\Gamma^{n}f_{\delta}(\Lambda_{\ast,\delta} (\psi)(s+\tau)) ds \\
   &= T^c(t) \big( \Pi^c (\Lambda_{\ast,\delta}(\psi)(\tau))\big)
      + \lim_{n\to \infty} \int_{0}^{t}T^c(t-s)\Pi^{c}\Gamma^{n}f_{\delta}(\tilde \Lambda(s)) ds \\
   &\quad -\lim_{n\to \infty} \int_{t}^{\infty}T^u(t-s)\Pi^{u}\Gamma^{n}f_{\delta}(\tilde\Lambda(s)) ds 
     + \lim_{n\to \infty} \int_{-\infty}^{t} T^s(t-s)\Pi^{s}\Gamma^{n}f_{\delta}(\tilde\Lambda(s)) ds \\
   &= {\cal F}_\delta \big( \Pi^c (\Lambda_{\ast,\delta}(\psi)(\tau)),\,\tilde\Lambda \big)(t), \quad t\in {\bb R},     
\end{align*}}%
that is, $\tilde\Lambda$ is a fixed point of 
${\cal F}_\delta \big(\Pi^c(\Lambda_{\ast,\delta}(\psi)(\tau)),\,\cdot \big)$. 
The uniqueness of the fixed points yields 
$ \tilde\Lambda = \Lambda_{\ast,\delta }\big( \Pi^c (\Lambda_{\ast,\delta}(\psi)(\tau)) \big)$, 
and hence
$$
   \Lambda_{\ast,\delta }(\psi )(t+\tau)=\tilde\Lambda(t)
   = \Lambda_{\ast,\delta }\big( \Pi^c (\Lambda_{\ast,\delta}(\psi)(\tau)) \big)(t), \quad t\in {\bb R},
$$
as desired. \qed

\vskip 5mm

For $\delta\in (0,\delta_1]$ let $F_{\ast,\delta}:E^{c}\to E^{su}$ be 
the map defined by $F_{\ast,\delta}(\psi):=\Pi^{su}\circ \mathrm{ev}_{0}\circ \Lambda_{\ast,\delta}(\psi)$ for $\psi\in E^{c}$, where $\mathrm{ev}_0$
is the evaluation map: $\mathrm{ev}_0(y):=y(0)$ for $y\in C({\bb R};X)$. 
Since  
\begin{align*}
   \Lambda_{\ast,\delta} (\psi)(0) 
   &= \psi-\lim_{n\to \infty} \int_{0}^{\infty}
       T^u(-s)\Pi^{u}\Gamma^{n}f_{\delta}(\Lambda_{\ast,\delta} (\psi)(s)) ds \\
   &\quad + \lim_{n\to \infty} \int_{-\infty}^{0}T^s(-s)\Pi^{s}\Gamma^{n} 
                 f_{\delta}(\Lambda_{\ast,\delta} (\psi)(s)) ds , 
\end{align*} 
it follows that 
\begin{align}\label{370}
   F_{\ast,\delta}(\psi)
   &= -\lim_{n\to \infty} \int_{0}^{\infty}
      T^u(-s)\Pi^{u}\Gamma^{n}f_{\delta}(\Lambda_{\ast,\delta}(\psi)(s))ds \notag \\
   &\quad + \lim_{n\to\infty}\int_{-\infty}^{0}
       T^s(-s)\Pi^{s}\Gamma^{n} f_{\delta}(\Lambda_{\ast,\delta} (\psi)(s))ds, 
       \quad \psi\in E^{c}; 
\end{align} 
and in particular 
$\Lambda_{\ast,\delta } (\psi)(0)=\psi+F_{\ast,\delta }(\psi)$ for 
$\psi\in E^{c}$. 

Let us set 
$$
   W^{c}_\delta
   :=\mathrm{graph}\,F_{\ast,\delta}=\big\{\psi+F_{\ast,\delta }(\psi): \psi\in E^c \big\}.
$$ 

\begin{Prop}\label{Proposition D3} 
The map $F_{\ast,\delta}$ and its graph $W^c_\delta$ have the 
following properties:
\begin{enumerate}
   \item \ $F_{\ast,\delta}$ is (globally) Lipschitz continuous, i.e.,
$$
    \|F_{\ast,\delta}(\psi_1) -F_{\ast,\delta}(\psi_2)\|_X
    \leq L(\delta)\|\psi_1-\psi_2\|_X, \quad \psi_1,\psi_2\in E^{c}, 
$$
         where $L(\delta):=4C^2C_1\zeta_\ast(\delta)/(\alpha -\eta)$. 
   \item \ Let $\hat\phi \in W^c_\delta$ and $\tau\in {\bb R}$. 
         Then the solution of ($E_\delta$) through $(\tau,\hat\phi)$, 
         $x(t;\tau,\hat\phi,f_\delta)$, exists on ${\bb R}$ and 
$$
    x_t(\tau,\hat\phi,f_\delta)=\Lambda _{\ast,\delta}(\hat\psi)(t-\tau), 
    \quad t\in {\bb R},  
$$
         where $\hat\psi=\Pi^c\hat\phi$. 
   \item \ Moreover for $\hat\phi \in W^c_\delta$ and $\tau\in {\bb R}$,
$$
   \Pi^{su} x_t(\tau,\hat\phi,f_\delta) 
   = F_{\ast,\delta}\big( \Pi^cx_t(\tau,\hat\phi,f_\delta) \big),
   \quad t\in{\bb R}. 
$$
In particular $W^c_\delta$ is invariant for ($E_\delta$), that is, 
$x_t(\tau,\hat\phi,f_\delta) \in W^c_\delta$ for $t\in {\bb R}$, 
provided that $\hat \phi\in W^c_\delta$.
\end{enumerate}
\end{Prop}

\noindent
{\it Proof.} (i) By (\ref{370}) and Proposition \ref{Proposition D2} (i)
\begin{align*}
   \big\|\Pi^s\big( F_{\ast,\delta}(\psi_1)-F_{\ast,\delta}(\psi_2) \big)\big\|_X
   &\leq \int_{-\infty}^0 CC_1e^{\alpha s}\zeta _\ast(\delta )
         \| \Lambda_{\ast,\delta } (\psi_1)(s)- \Lambda_{\ast,\delta } (\psi_2)(s) \|_X ds \\
   &\leq \int_{-\infty}^0 CC_1e^{\alpha s}\zeta _\ast(\delta )
         \| \Lambda_{\ast,\delta } (\psi_1)-\Lambda_{\ast,\delta }(\psi_2)\|_{Y_\eta}e^{\eta|s|} ds \\
   &= \frac{CC_1\zeta _\ast(\delta )}{\alpha -\eta} 
           \| \Lambda_{\ast,\delta } (\psi_1)-\Lambda_{\ast,\delta }(\psi_2)\|_{Y_\eta} \\
   &\leq  \frac{L(\delta)}{2}\|\psi_1 -\psi_2 \|_X . 
\end{align*}
A similar calculation gives
$$
   \big\|\Pi^u\big( F_{\ast,\delta}(\psi_1)-F_{\ast,\delta}(\psi_2) \big)\big\|_X
   \leq  \frac{L(\delta)}{2}\|\psi_1 -\psi_2 \|_X,  
$$
and therefore 
$$
  \|F_{\ast,\delta}(\psi_1)-F_{\ast,\delta}(\psi_2)\|_X \leq L(\delta)\|\psi_1 -\psi_2 \|_X.  
$$

\noindent
(ii) Applying Lemma \ref{Lemma D1} (i), we deduce that 
$\Lambda _{\ast,\delta}(\hat\psi) \in C({\bb R};X_0)$ 
and that the $X$-valued function 
$\xi(t):=\big(\Lambda _{\ast,\delta}(\hat\psi)(t)\big)[0]$ $(t\in {\bb R})$ 
satisfies $\xi_t=\Lambda _{\ast,\delta}(\hat\psi)(t)$ for $t\in {\bb R}$ 
and is a solution of $(E_\delta)$ on ${\bb R}$ with 
$\xi_0=\Lambda _{\ast,\delta}(\hat\psi)(0)=\hat\psi+F_{\ast,\delta }(\hat\psi)=\hat\phi$. Let $x(t):=\xi(t-\tau)$. 
Then $x(t)$ is a solution of $(E_\delta)$ on ${\bb R}$ with $x_\tau=\hat\phi$, 
so that $x(t)=x(t;\tau,\hat\phi,f_\delta)$ 
for $t\in {\bb R}$. Consequently, 
$$
   x_t(\tau,\hat\phi,f_\delta)=\xi_{t-\tau}
   =\Lambda _{\ast,\delta}(\hat\psi)(t-\tau), \qquad t\in {\bb R}.
$$

\noindent
(iii) Notice from Proposition \ref{Proposition D2} (ii) that 
$\Lambda_{\ast,\delta}(\hat\psi)(t-\tau)= \Lambda_{\ast,\delta }\big( \Pi^c (\Lambda_{\ast,\delta}(\hat\psi)(t-\tau)) \big)(0)$ 
for $\hat\psi:=\Pi^c\hat\phi$, which, combined with (ii),  
yields 
\begin{align*}
   \Pi^{su}x_t(\tau,\hat\phi,f_\delta)
   &= \Pi^{su}\big(\Lambda_{\ast,\delta }\big( \Pi^c (\Lambda_{\ast,\delta}(\hat\psi)(t-\tau)) \big)(0)\big) \\
   &= \Pi^{su}\big(\Lambda_{\ast,\delta }\big( \Pi^c x_t(\tau,\hat\phi,f_\delta)\big)(0) \big) \\
   &= F_{\ast,\delta}\big( \Pi^c x_t(\tau,\hat\phi,f_\delta)\big).
\end{align*}
The latter part of (iii) is obivous. \qed

\vskip 5mm

Now assume that $\Sigma^u=\emptyset$, i.e., $E^u=\{0\}$. 
Fix a $\delta\in (0,\delta_1]$ and let
$$
   K:=CC_1\zeta _\ast(\delta), \qquad\mu:=K+\varepsilon. 
$$

\begin{Prop}\label{Proposition D4}
Let $x(t)$ be a solution of ($E_\delta$) on an interval $J:=[t_0,t_1]$. 
Given $\tau\in J$, put $\hat\phi:=\Pi^cx_\tau+F_{\ast,\delta}(\Pi^cx_\tau)$. 
Then the following inequalities hold:
\begin{enumerate}
   \item For $t_0\leq t\leq \tau$ 
$$
   \|\Pi^cx_t-\Pi^cx_t(\tau,\hat\phi,f_\delta)\|_X
   \leq K\int_t^\tau e^{\mu(s-t)} \|\Pi^sx_s-\Pi^sx_s(\tau,\hat\phi,f_\delta)\|_Xds.
$$
   \item Moreover for $t_0\leq t\leq \tau$ 
$$
   \|\Pi^cx_t-\Pi^cx_t(\tau,\hat\phi,f_\delta)\|_X
   \leq K\int_t^\tau e^{\mu'(s-t)} \|\xi(s)\|_X ds,
$$
         where $\mu':=\mu+KL(\delta)$ and 
         $\xi(t):=\Pi^s x_t-F_{\ast,\delta}(\Pi^cx_t)$ for $t\in {\bb R}$.
\end{enumerate}
\end{Prop}

For the proof of Proposition \ref{Proposition D4}, we need the following 
lemma. 

\begin{Lem}\label{Lemma D2} Let $g(t)$, $h(t)$ and $r(t)$ are real 
valued continuous functions on the interval $[t_0,\tau]$ such that 
$r(t)\geq 0$ and  
\begin{align}\label{380}
   g(t)\leq \int_t^\tau h(s)ds + \int_t^\tau r(s)g(s)ds
\end{align}
for $t\in [t_0,\tau]$. Then we have  
$$
   g(t) \leq \int_t^\tau h(s) \exp\left( {\int_t^s r(\sigma)d\sigma}\right)ds, 
   \quad t\in [t_0,\tau].
$$
\end{Lem}

\noindent
{\it Proof.} Put 
$$
   F(t):=\int_t^\tau r(s)g(s)ds, \quad H(t):= \int_t^\tau h(s)ds
   \quad \text{and} \quad R(t):=\int_t^\tau r(s)ds.
$$
By the assumption (\ref{380}), $-F'(t)=r(t)g(t)\leq r(t)(F(t)+H(t))$ 
and hence 
$$
   \frac{d}{dt}\big( e^{-R(t)}F(t)\big)
   = e^{-R(t)}\big(r(t)F(t)+F'(t) \big)
   \geq -r(t)e^{-R(t)}H(t).
$$
Since $F(\tau)=H(\tau)=0$, 
\begin{align*}
   e^{-R(t)}F(t)
   &\leq \int_t^\tau r(s)e^{-R(s)}H(s)ds \\
   &= -e^{-R(t)}H(t) - \int_t^\tau e^{-R(s)}H'(s)ds \\
   &= -e^{-R(t)}H(t) + \int_t^\tau e^{-R(s)}h(s)ds,
\end{align*}
and hence
$$
   F(t)+H(t)\leq \int_t^\tau e^{R(t)-R(s)}h(s)ds.
$$
So, by using (\ref{380}) again, we get 
$$
   g(t) 
   \leq F(t)+H(t) 
   \leq \int_t^\tau h(s) \exp\left( {\int_t^s r(\sigma)d\sigma}\right)ds.
$$
This completes the proof. \qed

\vskip 5mm

\noindent
{\it Proof of Proposition \ref{Proposition D4}.} 
(i) By Proposition \ref{Proposition D3} (ii) and (iii), the solution 
$x(t;\tau,\hat\phi,f_\delta)$ exists on ${\bb R}$ and 
$\Pi^sx_t(\tau,\hat\phi,f_\delta)=F_{\ast,\delta}\big(\Pi^cx_t(\tau,\hat\phi,f_\delta)\big)$ for 
$t\in {\bb R}$. Since $t\leq \tau$, VCF gives 
$$
   x_\tau(\tau,\hat\phi,f_\delta)
   = T(\tau-t)x_t(\tau,\hat\phi,f_\delta)
     + \lim_{n\to\infty}\int_t^\tau T(\tau-s)\Gamma ^nf_\delta(x_s(\tau,\hat\phi,f_\delta))ds,
$$
in particular
$$
   \Pi^cx_\tau(\tau,\hat\phi,f_\delta)
   = T^c(\tau-t)\Pi^c x_t(\tau,\hat\phi,f_\delta)
     + \lim_{n\to\infty}\int_t^\tau T^c(\tau-s)\Pi^c\Gamma ^nf_\delta(x_s(\tau,\hat\phi,f_\delta))ds.
$$
By the group property of $\{T^c(t)\}_{t\in {\bb R}}$ 
\begin{align}\label{385}
   \Pi^cx_t(\tau,\hat\phi,f_\delta)
   &= T^c(t-\tau)\Pi^c x_\tau(\tau,\hat\phi,f_\delta) \notag \\
   &\quad - \lim_{n\to\infty}\int_t^\tau T^c(t-s)\Pi^c\Gamma ^nf_\delta(x_s(\tau,\hat\phi,f_\delta))ds.
\end{align}
Similarly for the solution $x(t)$ 
$$
   \Pi^cx_t= T^c(t-\tau)\Pi^c x_\tau
   - \lim_{n\to\infty}\int_t^\tau T^c(t-s)\Pi^c\Gamma ^nf_\delta(x_s)ds.
$$ 
Noting that $\Pi^c x_\tau(\tau,\hat\phi,f_\delta)=\Pi^c \hat\phi=\Pi^c x_\tau$, we obtain
\begin{align*}
   \|\Pi^cx_t - \Pi^cx_t(\tau,\hat\phi,f_\delta) \|_X
   &\leq \int_t^\tau CC_1e^{\varepsilon |t-s|}\zeta_\ast(\delta )\| x_s - x_s(\tau,\hat\phi,f_\delta)\|_Xds \\
   &\leq \int_t^\tau K e^{\varepsilon (s-t)}\big(\| \Pi^sx_s - \Pi^s x_s(\tau,\hat\phi,f_\delta)\|_X \\
   &\quad +\| \Pi^cx_s - \Pi^c x_s(\tau,\hat\phi,f_\delta)\|_X \big)ds, 
\end{align*}
so that
\begin{align*}
   e^{\varepsilon t}\|\Pi^cx_t - \Pi^cx_t(\tau,\hat\phi,f_\delta) \|_X
   &\leq \int_t^\tau Ke^{\varepsilon s} \| \Pi^sx_s - \Pi^s x_s(\tau,\hat\phi,f_\delta)\|_X ds \\
   &\quad +\int_t^\tau Ke^{\varepsilon s}\| \Pi^cx_s - \Pi^c x_s(\tau,\hat\phi,f_\delta)\|_X ds
\end{align*}
for $t_0\leq t\leq \tau$. Applying Lemma \ref{Lemma D2}, we get
$$
   e^{\varepsilon t} \|\Pi^cx_t - \Pi^cx_t(\tau,\hat\phi,f_\delta) \|_X
   \leq \int_t^\tau Ke^{K(s-t)} e^{\varepsilon s}\|\Pi^sx_s - \Pi^sx_s(\tau,\hat\phi,f_\delta) \|_Xds,
$$
which implies (i).

\noindent
(ii) By virtue of Proposition \ref{Proposition D3} (iii) and (i)
\begin{align*}
   \|\Pi^sx_s - \Pi^s x_s(\tau,\hat\phi,f_\delta) \|_X
   &\leq \|\Pi^s x_s-F_{\ast,\delta}(\Pi^c x_s)\|_X \\
   &\quad +\big\|F_{\ast,\delta}(\Pi^c x_s)- F_{\ast,\delta}\big(\Pi^c x_s(\tau,\hat\phi,f_\delta)\big)\big\|_X \\
   &\leq  \|\xi(s)\|_X + L(\delta)\| \Pi^c x_s - \Pi^c x_s(\tau,\hat\phi,f_\delta)\|_X
\end{align*}
for $s\in J$. Hence it follows from (i) that 
\begin{align*}
   e^{\mu t}\|\Pi^cx_t - \Pi^cx_t(\tau,\hat\phi,f_\delta) \|_X
   &\leq K\int_t^{\tau}e^{\mu s}\|\Pi^sx_s-\Pi^sx_s(\tau,\hat{\phi},f_{\delta})\|_Xds \\
   &\leq \int_t^\tau Ke^{\mu s}\|\xi(s)\|_Xds  \\
   &\quad +\int_t^\tau KL(\delta)e^{\mu s}\| \Pi^c x_s - \Pi^c x_s(\tau,\hat\phi,f_\delta)\|_Xds.
\end{align*}
Then another application of Lemma \ref{Lemma D2} readily yields (ii). \qed 

\vskip 5mm

Recall that 
\begin{align}\label{386}
   K:=CC_1\zeta_\ast(\delta), \quad 
   \mu:=K+\varepsilon, \quad 
   \mu':=\mu+KL(\delta)=K(1+L(\delta))+\varepsilon .
\end{align}

\begin{Prop}\label{Proposition D5}
Assume that $\Sigma^u=\emptyset$ and $x(t)$ is a solution of ($E_\delta$) 
on $J=[t_0,t_1]$. Define $\hat x_t\in W^c_\delta$ by 
$\hat x_t:=\Pi^c x_t+ F_{\ast,\delta}(\Pi^c x_t)$ for $t\in J$, and set 
$y(s;t):=\Pi^c x_s(t,\hat x_t,f_\delta)$ for $t\in J$ and $s\leq t$. 
Then the following inequality holds:
$$
   \|y(s;t)-y(s;t_0)\|_X
   \leq K\int_{t_0}^t e^{\mu'(\theta-s)}\|\xi(\theta)\|_Xd\theta, 
        \quad s\leq t_0,
$$
where $\xi(\theta):=\Pi^s x_\theta- F_{\ast,\delta}(\Pi^c x_\theta)$ 
for $\theta \in [t_0,t]$.
\end{Prop}

\noindent
{\it Proof.} Suppose that $s\leq t_0$. By the same reasoning as (\ref{385}) 
\begin{align}\label{387a}
   \Pi^cx_s(t,\hat x_t,f_\delta)
   = T^c(s-t)\Pi^c\hat x_t-\lim_{n\to\infty}\int_s^t 
     T^c(s-\sigma)\Pi^c\Gamma ^nf_\delta(x_\sigma(t,\hat x_t,f_\delta))d\sigma.
     \quad 
\end{align}
Applying VCF to $x_t$ and using $\Pi^c\hat x_\tau=\Pi^c x_\tau$ 
($\tau\in J$), we deduce that
$$
   \Pi^c \hat x_t=T^c(t-t_0)\Pi^c \hat x_{t_0}+\lim_{n\to\infty} 
   \int_{t_0}^t T^c(t-\sigma)\Pi^c\Gamma ^nf_\delta(x_\sigma)d\sigma, 
$$
and, thus, (\ref{387a}) becomes 
\begin{align*}
   \Pi^cx_s(t,\hat x_t,f_\delta)
   &= T^c(s-t_0)\Pi^c\hat x_{t_0} + \lim_{n\to\infty}\int_{t_0}^t 
      T^c(s-\sigma)\Pi^c\Gamma ^nf_\delta(x_\sigma)d\sigma \\
   &\quad -\lim_{n\to\infty}\int_s^{t} 
      T^c(s-\sigma)\Pi^c\Gamma ^nf_\delta(x_\sigma(t,\hat x_t,f_\delta))d\sigma, \quad t\in J. 
\end{align*}
Therefore
\begin{align}\label{387} 
   \|y(s;t)-y(s;t_0)\|_X 
   &= \|\Pi^c x_s(t,\hat x_t,f_\delta)- \Pi^cx_s(t_0,\hat x_{t_0},f_\delta) \|_X \notag \\
   &= \Big\| \lim_{n\to\infty}\int_{t_0}^t T^c(s-\sigma)\Pi^c\Gamma ^nf_\delta(x_\sigma)d\sigma \notag \\
   &\quad -\lim_{n\to\infty}\int_s^{t} T^c(s-\sigma)\Pi^c\Gamma ^nf_\delta(x_\sigma(t,\hat x_t,f_\delta))d\sigma \notag \\
   &\quad +\lim_{n\to\infty}\int_s^{t_0} T^c(s-\sigma)\Pi^c\Gamma ^nf_\delta(x_\sigma(t_0,\hat x_{t_0},f_\delta))d\sigma\,\Big\|_X \notag \\
   &\leq \int_{t_0}^t CC_1 e^{\varepsilon |s-\sigma|} \zeta_\ast(\delta)
                \| x_\sigma-x_\sigma(t,\hat x_t,f_\delta)\|_X d\sigma \notag \\
   &\quad +\int_s^{t_0} CC_1 e^{\varepsilon |s-\sigma|} \zeta_\ast(\delta)
                \|x_\sigma(t_0,\hat x_{t_0},f_\delta)  -x_\sigma(t,\hat x_t,f_\delta)\|_X d\sigma. 
\end{align}
Observe that
\begin{align}\label{390}
   \|x_\sigma-x_\sigma(t,\hat x_t,f_\delta)\|_X 
   &\leq \|\Pi^s x_\sigma- \Pi^s x_\sigma(t,\hat x_t,f_\delta)\|_X 
         + \|\Pi^c x_\sigma- \Pi^c x_\sigma(t,\hat x_t,f_\delta)\|_X \notag \\
   &\leq \|\Pi^s x_\sigma- F_{\ast,\delta}(\Pi^cx_\sigma)\|_X 
         + \|F_{\ast,\delta}(\Pi^cx_\sigma) -F_{\ast,\delta}(\Pi^c x_\sigma(t,\hat x_t,f_\delta))\|_X \notag \\
   &\quad +\|\Pi^c x_\sigma- \Pi^c x_\sigma(t,\hat x_t,f_\delta)\|_X \notag \\
   &\leq \|\xi(\sigma)\|_X+ (1+L(\delta))  \|\Pi^c x_\sigma- \Pi^c x_\sigma(t,\hat x_t,f_\delta)\|_X, 
\end{align}
where we used Proposition \ref{Proposition D3} (i) and (iii). 
Note also that
\begin{align}\label{395}
   \|x_\sigma(t_0,\hat x_{t_0},f_\delta)-x_\sigma(t,\hat x_t,f_\delta)\|_X 
   &\leq \|\Pi^s x_\sigma(t_0,\hat x_{t_0},f_\delta)- \Pi^s x_\sigma(t,\hat x_t,f_\delta)\|_X  \notag \\
   &\quad +\|\Pi^c x_\sigma(t_0,\hat x_{t_0},f_\delta)- \Pi^c x_\sigma(t,\hat x_t,f_\delta)\|_X  \notag \\
   &= \| F_{\ast,\delta} (\Pi^c x_\sigma(t_0,\hat x_{t_0},f_\delta))
      - F_{\ast,\delta}(\Pi^cx_\sigma(t,\hat x_{t},f_\delta))\|_X  \notag \\
   &\quad +\|\Pi^c x_\sigma(t_0,\hat x_{t_0},f_\delta)- \Pi^c x_\sigma(t,\hat x_t,f_\delta)\|_X  \notag \\
   &\leq (1+L(\delta)) \|\Pi^c x_\sigma(t_0,\hat x_{t_0},f_\delta)- \Pi^c x_\sigma(t,\hat x_t,f_\delta)\|_X  \notag \\
   &= (1+L(\delta))\|y(\sigma ;t)-y(\sigma ;t_0)\|_X. 
\end{align}
In view of (\ref{387}), (\ref{390}) and (\ref{395}), combined with Proposition \ref{Proposition D4} (ii), we deduce
{\allowdisplaybreaks
\begin{align}\label{400} 
   \|y(s;t)-y(s;t_0)\|_X 
   &\leq \int_{t_0}^t K e^{\varepsilon (\sigma-s)} 
          \big( \|\xi(\sigma)\|_X+ (1+L(\delta))  \|\Pi^c x_\sigma- \Pi^c x_\sigma(t,\hat x_t,f_\delta)\|_X\big)d\sigma  \notag \\
   &\quad +\int_s^{t_0} K e^{\varepsilon (\sigma-s)} 
          (1+L(\delta)) \|y(\sigma ;t)-y(\sigma ;t_0)\|_X d\sigma  \notag \\
   &\leq  \int_{t_0}^t K e^{\varepsilon (\sigma-s)}  \|\xi(\sigma)\|_X d\sigma  \notag \\
   &\quad +\int_{t_0}^t K e^{\varepsilon (\sigma-s)} 
          (1+L(\delta))\Big( K \int_\sigma^t e^{\mu'(\tau-\sigma)}\|\xi(\tau)\|_Xd\tau \Big)d\sigma  \notag \\
   &\quad +\int_s^{t_0} K e^{\varepsilon (\sigma-s)} 
          (1+L(\delta)) \|y(\sigma ;t)-y(\sigma ;t_0)\|_X d\sigma. 
\end{align}}%
Notice that the second term of the right-hand side becomes 
$$
  K\int_{t_0}^t \big(e^{\varepsilon(t_0-s)+\mu'(\sigma-t_0)}-e^{\varepsilon (\sigma-s)} \big)\|\xi(\sigma)\|_Xd\sigma  
$$
because of (\ref{386}). So we see from (\ref{400}) that for $s\leq t_0$ 
\begin{align*}
   e^{\varepsilon s}\|y(s;t)-y(s;t_0)\|_X 
   &\leq K\int_{t_0}^t e^{(\varepsilon-\mu')t_0+\mu'\sigma}\|\xi(\sigma)\|_Xd\sigma  \\
   &\quad +K(1+L(\delta)) \int_s^{t_0} e^{\varepsilon \sigma} \|y(\sigma ;t)-y(\sigma ;t_0)\|_X d\sigma. 
\end{align*}
By Gronwall's inequality and (\ref{386})
\begin{align*}
   e^{\varepsilon s}\|y(s;t)-y(s;t_0)\|_X 
   &\leq  \Big( K\int_{t_0}^t e^{(\varepsilon-\mu')t_0+\mu'\sigma}\|\xi(\sigma)\|_Xd\sigma  \Big) e^{K(1+L(\delta))(t_0-s)}   \\
   &= K e^{-(\mu'-\varepsilon )s} \int_{t_0}^t e^{\mu'\sigma}\|\xi(\sigma)\|_Xd\sigma,       
\end{align*}
and therefore
$$
   \|y(s;t)-y(s;t_0)\|_X \leq K  \int_{t_0}^t e^{\mu'(\sigma-s)}\|\xi(\sigma)\|_Xd\sigma, \quad s\leq t_0, 
$$
as required. \qed

\begin{Prop}\label{Proposition D6} 
Assume that $\Sigma^u=\emptyset$ and $\delta \in (0,\delta_1]$ satisfies 
\begin{align}\label{405}
   \max\left( \mu', \frac{K(\alpha -\varepsilon)}{\alpha -\mu'}\right) 
   <\alpha. 
\end{align}
If $x(t)$ is a solution of ($E_\delta$) on $J=[t_0,t_1]$, then the 
function $\xi(t):=\Pi^s x_t-F_{\ast,\delta}(\Pi^c x_t)$ satisfies the 
inequality
$$
   \|\xi(t)\|_X\leq C\|\xi(t_0)\|_Xe^{-\beta_0(t-t_0)}, \quad t\in J,
$$  
where $\beta _0:=\alpha -K(\alpha- \varepsilon) /(\alpha -\mu')>0$. 
If in particular $J=[t_0,\infty)$, ${\rm dist}\,(x_t,W^c_\delta)$ 
tends to $0$ exponentially as $t\to\infty$. 
\end{Prop}

\noindent
{\it Proof.} Observe that for $t\in J$ 
{\allowdisplaybreaks
\begin{align*}
   \xi(t)-T^s(t-t_0)\xi(t_0) 
   &= \Pi^s x_t-F_{\ast,\delta}(\Pi^c x_t)-T^s(t-t_0)(\Pi^s x_{t_0}-F_{\ast,\delta}(\Pi^c x_{t_0})) \\
   &= \Pi^s (x_t-T(t-t_0)x_{t_0})-F_{\ast,\delta}(\Pi^c x_t)+ T^s(t-t_0)F_{\ast,\delta}(\Pi^c x_{t_0}) \\
   &= \lim_{n\to\infty}\int_{t_0}^t T^s(t-s)\Pi^s\Gamma ^nf_\delta(x_s)ds \qquad \quad (\,{\rm by \ VCF}\,) \\
   &\quad -\lim_{n\to\infty}\int_{-\infty}^0 T^s(-s)\Pi^s\Gamma ^nf_\delta( \Lambda_{\ast,\delta}(\Pi^c x_t)(s))ds \\
   &\quad +\lim_{n\to\infty}\int_{-\infty}^0 T^s(t-t_0-s)\Pi^s\Gamma ^nf_\delta(\Lambda_{\ast,\delta}(\Pi^c x_{t_0})(s))ds \\
   &= \lim_{n\to\infty}\int_{t_0-t}^0 T^s(-s)\Pi^s\Gamma ^nf_\delta(x_{s+t})ds  \\
   &\quad -\lim_{n\to\infty}\int_{-\infty}^0 T^s(-s)\Pi^s\Gamma ^nf_\delta( \Lambda_{\ast,\delta}(\Pi^c x_t)(s))ds \\
   &\quad +\lim_{n\to\infty}\int_{-\infty}^{t_0-t} T^s(-s)\Pi^s\Gamma ^nf_\delta(\Lambda_{\ast,\delta}(\Pi^c x_{t_0})(t-t_0+s))ds \\
   &= \lim_{n\to\infty}\int_{t_0-t}^0 T^s(-s)\Pi^s\Gamma ^n\big(f_\delta(x_{s+t})- f_\delta( \Lambda_{\ast,\delta}(\Pi^c x_t)(s)) \big)ds \\
   &\quad +\lim_{n\to\infty}\int_{-\infty}^{t_0-t} T^s(-s)\Pi^s\Gamma ^n
   \big( f_\delta(\Lambda_{\ast,\delta}(\Pi^c x_{t_0})(t-t_0+s)) \\
   &\quad  -f_\delta( \Lambda_{\ast,\delta}(\Pi^c x_t)(s)) \big)ds.
\end{align*}}%
If we set $\hat x_t:=\Pi^c x_t+F_{\ast,\delta}(\Pi^c x_t)$ for 
$t\in J$, by Proposition \ref{Proposition D3} (ii)   
$\Lambda_{\ast,\delta}(\Pi^c x_t)(s)=x_s(0,\hat x_t,f_\delta)=x_{s+t}(t,\hat x_t,f_\delta)$ and 
$\Lambda_{\ast,\delta}(\Pi^c x_{t_0})(t-t_0+s)=x_{t-t_0+s}(0,\hat x_{t_0},f_\delta)=x_{s+t}(t_0,\hat x_{t_0},f_\delta)$ 
in particular for $s\in {\bb R}^-$. So 
\begin{align*}
   \xi(t) 
   &= T^s(t-t_0)\xi(t_0) 
      +\lim_{n\to\infty}\int_{t_0-t}^0 T^s(-s)\Pi^s\Gamma ^n 
      \big(f_\delta(x_{s+t})- f_\delta(x_{s+t}(t,\hat x_t,f_\delta)) \big)ds \\
   &\quad +\lim_{n\to\infty}\int_{-\infty}^{t_0-t} T^s(-s)\Pi^s\Gamma ^n
          \big( f_\delta(x_{s+t}(t_0,\hat x_{t_0},f_\delta))- f_\delta(x_{s+t}(t,\hat x_t,f_\delta)) \big)ds,
\end{align*}
and thus
\begin{align*}
   \|\xi(t)\|_X 
   &\leq Ce^{-\alpha (t-t_0)}\|\xi(t_0)\|_X 
         + \int_{t_0-t}^0 Ke^{\alpha s}\|x_{s+t}-x_{s+t}(t,\hat x_t,f_\delta) \|_Xds  \\
   &\quad  +\int_{-\infty}^{t_0-t} Ke^{\alpha s} 
         \|x_{s+t}(t_0,\hat x_{t_0},f_\delta)- x_{s+t}(t,\hat x_t,f_\delta)\|_X ds  \\
   &= Ce^{-\alpha (t-t_0)}\|\xi(t_0)\|_X 
         + \int_{t_0}^t Ke^{\alpha (\theta-t)}
           \|x_{\theta}-x_{\theta}(t,\hat x_t,f_\delta) \|_Xd\theta  \\
   &\quad +\int_{-\infty}^{t_0} Ke^{\alpha (\theta-t)} 
           \|x_{\theta}(t_0,\hat x_{t_0},f_\delta)- x_{\theta}(t,\hat x_t,f_\delta)\|_X d\theta. 
\end{align*} 
Since $x_{\theta}(t,\hat x_t,f_\delta)$ $(t\in J)$ can be written as 
\begin{align*}
   x_{\theta}(t,\hat x_t,f_\delta)
   &= \Pi^c x_{\theta}(t,\hat x_t,f_\delta) 
      + \Pi^s x_{\theta}(t,\hat x_t,f_\delta) \\
   &= \Pi^c x_{\theta}(t,\hat x_t,f_\delta) 
      + F_{\ast,\delta}( \Pi^c x_{\theta}(t,\hat x_t,f_\delta)),
      \quad \theta \in {\bb R}
\end{align*} 
(Proposition \ref{Proposition D3} (iii)), it follows from 
Proposition \ref{Proposition D3} (i) and 
Proposition \ref{Proposition D5} that for $\theta \leq t_0$
\begin{align*}
   \|x_{\theta}(t_0,\hat x_{t_0},f_\delta)-x_{\theta}(t,\hat x_t,f_\delta)\|_X 
   &\leq \|\Pi^cx_{\theta}(t_0,\hat x_{t_0},f_\delta)- \Pi^cx_{\theta}(t,\hat x_t,f_\delta)\|_X  \\
   &\quad + \|F_{\ast,\delta}(\Pi^c x_{\theta}(t_0,\hat x_{t_0},f_\delta))
          - F_{\ast,\delta}(\Pi^c x_{\theta}(t,\hat x_t,f_\delta)) \|_X   \\
   &\leq (1+L(\delta)) \|\Pi^c x_{\theta}(t_0,\hat x_{t_0},f_\delta)- \Pi^c x_{\theta}(t,\hat x_t,f_\delta)\|_X \\
   &= (1+L(\delta)) \|y(\theta;t)-y(\theta;t_0)\|_X \\
   &\leq (1+L(\delta)) K\int_{t_0}^t e^{\mu'(\tau-\theta)}\|\xi(\tau)\|_Xd\tau,
\end{align*} 
where $y(\theta;t)$ $(t\in J)$ is the one in 
Proposition \ref{Proposition D5}. On the other hand, for 
$t_0\leq \theta \leq t$   
\begin{align*}
   \|x_{\theta}-x_{\theta}(t,\hat x_t,f_\delta) \|_X
   &\leq \|\Pi^s x_{\theta}-\Pi^s x_{\theta}(t,\hat x_t,f_\delta) \|_X
          + \|\Pi^c x_{\theta}- \Pi^c x_{\theta}(t,\hat x_t,f_\delta) \|_X  \\
   &\leq \|\Pi^s x_{\theta}- F_{\ast,\delta}(\Pi^c x_{\theta}) \|_X
          +\|F_{\ast,\delta}(\Pi^c x_{\theta})- F_{\ast,\delta}(\Pi^c x_{\theta}(t,\hat x_t,f_\delta)) \|_X \\
   &\quad +\|\Pi^c x_{\theta}- \Pi^c x_{\theta}(t,\hat x_t,f_\delta) \|_X  \\  
   &\leq \|\xi(\theta)\|_X +(1+L(\delta)) \|\Pi^c x_{\theta}- \Pi^c x_{\theta}(t,\hat x_t,f_\delta) \|_X \\
   &\leq \|\xi(\theta)\|_X +(1+L(\delta)) K\int_{\theta}^t e^{\mu'(\sigma-\theta)}\|\xi(\sigma)\|_Xd\sigma,       
\end{align*}
where we used Proposition \ref{Proposition D3} (i), (iii) and 
Proposition \ref{Proposition D4} (ii). Thus we have 
\begin{align*} 
   \|\xi(t)\|_X 
   &\leq  Ce^{-\alpha (t-t_0)}\|\xi(t_0)\|_X  \\
   &\quad + \int_{t_0}^t Ke^{\alpha (\theta-t)}
          \Big( \|\xi(\theta)\|_X 
          +(1+L(\delta)) K\int_{\theta}^t e^{\mu'(\sigma-\theta)}\|\xi(\sigma)\|_Xd\sigma\Big) d\theta  \\ 
   &\quad + \int_{-\infty}^{t_0} Ke^{\alpha (\theta-t)} 
          (1+L(\delta)) K\Big(\,\int_{t_0}^t e^{\mu'(\tau-\theta)}\|\xi(\tau)\|_Xd\tau \Big)d\theta   \\
   &= Ce^{-\alpha (t-t_0)}\|\xi(t_0)\|_X 
        +\Big(K+\frac{K^2(1+L(\delta))}{\alpha -\mu'}\Big) \int_{t_0}^t e^{\alpha (\sigma-t)}\|\xi(\sigma)\|_X d\sigma,     
\end{align*} 
so that 
$$
   e^{\alpha t} \|\xi(t)\|_X 
   \leq Ce^{\alpha t_0} \|\xi(t_0)\|_X 
        +\hat K  \int_{t_0}^t e^{\alpha \sigma}\|\xi(\sigma)\|_X d\sigma,
$$
where $\hat K:=K+K^2(1+L(\delta))/(\alpha -\mu')$. An application of 
Gronwall's inequality gives 
$$
   e^{\alpha t} \|\xi(t)\|_X
   \leq Ce^{\alpha t_0} \|\xi(t_0)\|_X e^{\hat K(t-t_0)},
$$
and hence 
$$
   \|\xi(t)\|_X 
   \leq C \|\xi(t_0)\|_X e^{-(\alpha -\hat K)(t-t_0)}, \quad t\in J,
$$
which is the desired one because in view of (\ref{386}) 
$$
   \hat K=K\frac{\alpha -\varepsilon }{\alpha -\mu'}=\alpha -\beta_0. 
$$

The latter part of the proposition is evident. This completes the proof. \qed

\vskip 5mm
\noindent
{\it Proof of Theorem \ref{Theorem D1}.} The properties (ii) and (iii) of 
Theorem \ref{Theorem D1} are now immediate consequences of 
Propositions \ref{Proposition D3} and \ref{Proposition D6}, 
respectively. We verify the property (i). Observe that $Y_\eta$ is a 
subspace of $Y_{\eta'}$ if $\eta<\eta'<\alpha$, and denote the 
inclusion map by ${\cal J}:Y_{\eta}\to Y_{\eta'}$. It will be shown, 
in the appendix (Proposition \ref{Proposition D10}), that 
${\cal J}\Lambda_{\ast,\delta}$ is $C^1$ smooth as a map from $E^c$ 
to $Y_{\eta'}$; and hence 
$F_{\ast,\delta}=\Pi^{su}\circ \mathrm{ev}_{0}\circ {\cal J}\Lambda_{\ast,\delta}$ is also  $C^1$ smooth. Moreover, the relation
$$
   \big[[D({\cal J}\Lambda_{\ast,\delta})(0)](t)\big]\psi
   =T^c(t)\psi, \quad \psi\in E^c,\ t\in {\bb R}
$$
holds since $Df_\delta(0)=Df(0)=0$ (cf.\,(\ref{553}) and (\ref{553b})). 
In particular
$$
   DF_{\ast,\delta}(0)\psi
   = D(\Pi^{su}\circ \mathrm{ev}_{0}\circ {\cal J}\Lambda_{\ast,\delta})(0)\psi
   = \Pi^{su}T^c(0)\psi=\Pi^{su}\psi=0, \qquad \psi\in E^c,
$$
so that $ DF_{\ast,\delta}(0)=0$, which implies (i).  \qed

\vskip 5mm
\subsection{Stability for integral equations via the central equation}
In this subsection, introducing some ordinary differential equation 
which we call the central equation, we will study stability properties 
for the zero solution of Eq. ($E$). 

Assume that $\Sigma ^c \not=\emptyset$. Let $\{\phi_{1},\ldots,\phi_{d_c}\}$ 
be a basis for $E^{c}$, where $d_c$ is the dimension of $E^{c}$. 
Then based on the formal adjoint theory for Eq. (\ref{126c}) 
developed in \cite{murmat} (also refer to \cite{murmatnagnak})
where the formal adjoint theory is accomplished for Volterra 
difference equations), one can consider its dual basis as elements 
in the Banach space 
$$
   X^\sharp :=L^1_\rho({\bb R}^+; ({\bb C}^\ast)^m )
   = \big\{\psi:{\bb R}^+\to ({\bb C}^\ast)^m: 
     \text{$\psi(\tau)e^{-\rho \tau}$ is integrable on ${\bb R}^+$} \big\}
$$
with norm
$$
   \|\psi\|_{X^\sharp}
   := \int_0^\infty |\psi(\tau)|e^{-\rho \tau}d\tau, \quad \psi\in X^\sharp,
$$
where $({\bb C}^\ast)^m$ is the space of  $m$-dimensional row vectors 
with complex components equipped with the norm which is compatible 
with the one in ${\bb C}^m$, that is, $|z^\ast z|\leq |z^\ast|\,|z|$ for 
$z^\ast \in ({\bb C}^\ast)^m$ and $z\in {\bb C}^m$. To be more precise, 
if we set 
$$
   \lp \psi,\phi\rp
   := \int_{-\infty}^0 \left( \int_\theta^0 \psi(\xi-\theta)
      K(-\theta)\phi(\xi)d\xi \right)d\theta, \quad 
      (\psi,\phi) \in X^\sharp\times X,
$$ 
then this pairing defines a bounded bilinear form on $ X^\sharp\times X$ 
with the property 
$$
   |\lp \psi,\phi\rp| 
   \leq \|K\|_{\infty,\rho } \|\psi\|_{X^\sharp} \|\phi\|_X, 
        \quad (\psi,\phi) \in X^\sharp\times X; 
$$
here we recall that 
$ \|K\|_{\infty,\rho}=\mathrm{ess}\,\sup \{\|K(t)\|e^{\rho t} : t\geq 0 \}$. 
Then there exist $\{\psi_{1},\ldots,\psi_{d_c}\}$, elements of 
$X^\sharp$, such that $\lp \psi_{i}, \phi_{j} \rp=1$ if $i=j$ and 
$0$ otherwise, and $\lp \psi_{i},\phi\rp= 0$ for $\phi \in E^s$ and 
$i=1,2,\ldots,d_c$; we call  $\{\psi_{1},\ldots,\psi_{d_c}\}$ the dual 
basis of $\{\phi_{1},\ldots,\phi_{d_c}\}$. (See \cite{murmat} for details.)
 
Denote by $\Phi_c$ and $\Psi_c$, $(\phi_{1},\ldots,\phi_{d_c})$ and 
$\,^t(\psi_{1},\ldots,\psi_{d_c})$, the transpose of 
$(\psi_{1},\ldots,\psi_{d_c})$, respectively. Then, for any $\phi\in X$ 
the coordinate of its $E^{c}$-component with respest to the basis 
$\{\phi_{1},\ldots,\phi_{d_c}\}$, or $\Phi_c$ for short, is given by 
$$
   \lp \Psi_c,\phi\rp
   := {^t}\big(\lp\psi_1,\phi\rp, \ldots, \lp\psi_{d_c},\phi \rp\big) 
      \in {\bb C}^{d_c}, 
$$
and therefore the projection $\Pi^{c}$ is expressed, in terms of 
the basis $\Phi_c$ and its dual basis $\Psi_c$, by
\begin{align}\label{330a}
   \Pi^{c}\phi=\Phi_c\lp \Psi_c,\phi \rp, \quad \phi\in X.
\end{align}

Since $\{T^{c}(t)\}_{t\geq 0}$ is a strongly continuous semigroup on 
the finite dimensional space $E^{c}$, there exists a $d_c\times d_c$ 
matrix $G_c$ such that 
\begin{align}\label{330b}
   T^{c}(t)\Phi_c=\Phi_c e^{tG_c}, \quad t\geq 0,
\end{align}
and $\sigma(G_c)$, the spectrum of $G_c$, is identical with $\Sigma^{c}$. 
The $E^{c}$-components of solutions of Eq.($E_\delta $) can be described 
by a certain ordinary differential equation in ${\bb C}^{d_c}$. 
More precisely, let $x(t)$ be a solution of Eq.($E_\delta $) through 
$(\sigma,\phi)$, that is, $x(t)=x(t;\sigma,\phi,f)$. If we denote by 
$z_c(t)$ the component of  $\Pi^{c}x_t$ with respect to the basis $\Phi_c$, 
that is, $\Phi_c z_c(t):=\Pi^{c}x_t$, or $z_c(t):=\lp \Psi_c,x_t\rp$, then 
by virtue of \cite[Theorem 7]{mur} $z_c(t)$ satisfies the ordinary 
differential equation 
\begin{align}\label{331}
   \dot{z_c}(t)=G_cz_c(t)+ H_c f_\delta(\Phi_c z_c(t)+\Pi^{su}x_{t}), 
\end{align}
where $H_c$ is the $d_c\times m$ matrix such that  
$$
   H_c\,x:=\lim_{n\to\infty} \lp \Psi_c,\Gamma^n x \rp, \quad x\in {\bb C}^m.
$$  
Thus the $E^{c}$-components of solutions of Eq.($E_\delta$) are 
determined by solutions of (\ref{331}).  

In connection with Eq. (\ref{331}), let us consider the ordinary 
differential equations on ${\bb C}^{d_c}$
\begin{align}
   \dot{z}(t)=G_cz(t)+ H_c f_\delta(\Phi_c z(t)+F_{\ast,\delta}(\Phi_c z(t))) 
   \tag{$CE_\delta$}
\end{align}
and 
\begin{align}
   \dot{z}(t)=G_cz(t)+ H_c f(\Phi_c z(t)+F_{\ast}(\Phi_c z(t))). 
   \tag{$CE$}
\end{align}
We call Eq. ($CE$) (resp. Eq. ($CE_{\delta}$)) the central equation of 
($E$) (resp. ($E_{\delta}$)).

\begin{Prop}\label{Proposition D7} The following statements hold true:
\begin{enumerate}
   \item Let $x$ be a solution of Eq.($E_\delta$) on an interval $J$ 
         such that $x_{t}\in W^{c}_\delta$ $(t\in J)$. 
         Then the function $z_c(t):=\lp \Psi_c,x_t\rp$ satisfies the 
         equation ($CE_\delta$) on $J$. \\
         Conversely, if $z(t)$ satisfies the equation ($CE_\delta$) 
         on an interval $J$, then there exists a unique solution $x$ 
         of Eq.($E_\delta$) on $J$ such that $x_{t}\in W^{c}_\delta$ 
         and $\Pi^{c}x_{t}=\Phi_c z(t)$ on $J$.

   \item Let $x$ be a solution of Eq.($E$) on an interval $J$ such that 
         $x_{t}\in W^{c}_{\rm loc}(r,\delta)$ $(t\in J)$. Then the 
         function $z_c(t):=\lp \Psi_c,x_t\rp$ satisfies the equation 
         ($CE$) on $J$, together with the inequality $\sup_{t\in J}
         \|\Phi_cz_c(t)\|_X\leq r$. \\
         Conversely, if $z(t)$ satisfies the equation ($CE$) on an 
         interval $J$ together with the inequality $\sup_{t\in J}
         \|\Phi_cz(t)\|_X\leq r$, then there exists a unique solution $x$ 
         of Eq.($E$) on $J$ such that $x_{t}\in W^{c}_{\rm loc}(r,\delta)$ 
         and $\Pi^{c}x_{t}=\Phi_c z(t)$ on $J$.
\end{enumerate}
\end{Prop}

\noindent
{\it Proof.} 
Let us recall that whenever $\phi\in X$ satisfies $\|\Pi^c\phi\|_X\leq r$, 
$\phi\in W_{\delta}^c$ if and only if $\phi\in W^c_{\rm loc}(r,\delta)$; 
and consequently $F_{\ast,\delta}(\Pi^c\phi)=F_{\ast}(\Pi^c\phi)$ and 
$f_{\delta}(\Pi^c\phi+F_{\ast,\delta}(\Pi^c\phi))=
f(\Pi^c\phi+F_{\ast}(\Pi^c\phi))$. 
Therefore $z(t)$ satisfying $\sup_{t\in J}\|\Phi_cz(t)\|_X\leq r$ is 
a solution of ($CE_{\delta}$) on $J$ if and only if it is a solution of 
($CE$) on $J$.  Thus (ii) is a direct consequence of (i); 
so, in what follows, we will prove (i) only.

The former part of (i) directly follows from 
Proposition \ref{Proposition D3} (iii). Conversely, let $z(t)$ be a 
solution of ($CE_\delta$) on $J$. Pick $\tau\in J$ and set 
$\hat\phi= \Phi_c z(\tau)+F_{\ast,\delta}(\Phi_c z(\tau))$. 
Then it follows from Proposition \ref{Proposition D3} (iii) again that 
$x(t):=x(t;\tau,\hat\phi,f_\delta)$ is a solution of ($E_\delta$) on 
$J$ such that $x_t\in W^c_\delta$ for $t\in J$. By the former part, 
the function $z_c(t)$, defined by $\Phi_cz_c(t)=\Pi^c x_t$, is also 
a solution of $(CE_\delta)$ satisfying 
$z_c(\tau)=\lp \Psi_c, \hat\phi\rp= \lp \Psi_c, \Phi_c z(\tau)\rp =z(\tau)$. 
The uniqueness of solutions of $(CE_\delta)$ yields $z_c(t)=z(t)$ for 
$t\in J$ and hence $\Pi^c x_t=\Phi_cz_c(t)= \Phi_c z(t)$ for $t\in J$. 
\qed

\vskip 5mm

Since $f(0)=f_{\delta}(0)=0$, both equations ($CE$) and ($CE_\delta$) 
(as well as ($E$) and ($E_\delta$)) possess the zero solution.  
Notice that the zero solution of ($CE$) (resp. ($E$)) is uniformly 
asymptotically stable if and only if the zero solution of 
($CE_{\delta}$) (resp. ($E_\delta$)) is uniformly asymptotically stable. 
Likewise, the zero solution of ($CE$) (resp. ($E$)) is unstable 
if and only if the zero solution of ($CE_{\delta}$) (resp. ($E_\delta$)) 
is unstable. Here, for the definition of several stability properties 
utilized in this paper, we  refer readers to the books \cite{yos, hen}. 

Now suppose that $\Sigma^u=\emptyset$. Then the dynamics near the 
zero solution of ($E$) is determined by the dynamics near $z_c=0$ of 
($CE$) in the following sense.  

\begin{Thm}\label{Theorem D3} Assume that $\Sigma^u=\emptyset$. 
If the zero solution of ($CE$) is uniformly asymptotically stable 
(resp. unstable), then the zero solution of ($E$) is also 
uniformly asymptotically stable (resp. unstable). 
\end{Thm}

\noindent
{\it Proof.} By the fact stated in the preceding paragraph of the 
theorem, it is sufficient to establish that the uniform asymptotic 
stability (resp. instability) of the zero solution of ($CE_{\delta}$) 
implies the uniform asymptotic stability (resp. instability) of 
the zero solution of $(E_\delta)$. 

If the zero solution of ($CE_\delta$) is unstable, the instability 
of the zero solution of $(E_\delta)$ immediately follows from the 
invariance of $W^c_\delta$ (Proposition \ref{Proposition D3} (iii)). 
In what follows, under the assumption that the the zero solution of 
($CE_\delta$) is uniformly asymptotically stable, we will establish 
the uniform asymptotic stability of the zero solution of $(E_\delta)$, 
employing an idea utilized for the stability problems of parabolic 
partial differential equations in \cite[Theorem 6.1.4]{hen}.

By virtue of \cite[Theorem 4.2.1]{hen}, there exist positive 
constants $a$, $\bar K$ and a Liapunov function $V$ defined on 
$S_a:=\{y\in {\bb C}^{d_c}: |y|\leq a \}$ satisfying the following properties: 
\begin{enumerate}
   \item There exists a $b\in C({\bb R}^+;{\bb R}^+)$ which is strictly 
         increasing with $b(0)=0$ and 
$$
    b(|y|)\leq V(y)\leq |y| 
    \quad \text{for} \quad y\in S_a.
$$ 
   \item $|V(y)-V(z)| \leq {\bar K} |y-z|$ \ for \ $y,z\in S_a$.
   \item $\dot V(z)\leq -V(z)$ \ for \ $z\in S_a$, where $\dot V(z)$ 
         is defined by 
$$
   \dot V(z):=\limsup_{h\to +0}\frac{V(y(h))-V(z)}{h},
$$
         $y(h)$ being the solution of ($CE_\delta$) with $y(0)=z$.
\end{enumerate}   
Choose a positive number $\tau_0$ such that 
\begin{align}
   e^{-\tau_0}\leq \frac{1}{2} 
   \quad \text{and} \quad Ce^{-\beta_0 \tau_0} \leq \frac{1}{4},
\end{align} 
where $\beta_0$ is the one in Proposition \ref{Proposition D6}, 
and we may assume that $\beta_0>\mu'$, taking $\delta$ so small if 
necessary. Put $K_\infty:=\|K\|_{\infty,\rho}$ and take a positive 
number $P$ in such a way that 
\begin{align}
   P>\max\Big(1,\, \frac{4}{\beta_0-\mu'}\bar K KK_\infty \|\Psi_c\|\Big),
\end{align}  
and set
$$
   a_0:= \frac{ae^{-\eta \tau_0}}{4CK_\infty \|\Psi_c\|},
$$
where $\|\Psi_c\|:=\big(\sum_{j=1}^{d_c} \|\psi_j\|_{X^\sharp}^2\big)^{1/2}$. 
Let $\Omega$ be a neighborhood of $0$ in $X$ such that 
$$
   \lp \Psi_c,\phi\rp \in S_a, \quad \|\Pi^c \phi\|_X \leq a_0, 
   \quad \text{and} \quad Q\leq b(a)
$$
for $\phi\in \Omega$, where
$$
   Q:=V(\lp\Psi_c, \phi\rp)
      +\left( PC+\frac{\bar KK_\infty\|\Psi_c\|KC}{\beta_0-\mu'}\right)
      \big(\|\Pi^s\phi\|_X+\|F_{\ast,\delta}(\Pi^c\phi)\|_X\big),
$$ 
and consider the function $W(\phi)$ on $\Omega$ defined by 
$$
   W(\phi):=V(\lp \Psi_c,\phi\rp)+P\|\Pi^s\phi-F_{\ast,\delta}(\Pi^c\phi)\|_X,
   \quad \phi\in \Omega.
$$
$W$ is continuous in $\Omega$ with $W(0)=0$ and is positive in 
$\Omega\backslash \{0\}$ because of (i) and (ii).   

We will first certify the following claim.

\begin{Claim}
There exists a positive number $c_0$ such that, for any 
$t_0\in {\bb R}^+$ and $\phi\in X$ with $W(\phi)\leq c_0$, the solution 
$x(t;t_0,\phi,f_\delta)$ exists on $[t_0,t_0+\tau_0]$ and satisfies 
$x_t(t_0,\phi,f_\delta)\in \Omega$ for $t\in [t_0,t_0+\tau_0]$; 
in particular, $\|\Pi^c  x_t(t_0,\phi,f_\delta)\|_X\leq a_0$ in this 
interval.
\end{Claim}

Indeed, suppose that $x_t(t_0,\phi,f_\delta)$ is defined on the interval 
$[t_0,t_0+t_\ast )$ with $t_\ast \leq \tau_0$. Then    
$$
   x_t(t_0,\phi,f_\delta)
   = T(t-t_0)\phi+\lim_{n\to\infty}\int_{t_0}^t 
     T(t-s)\Gamma^n f_\delta(x_s(t_0,\phi,f_\delta))ds 
$$
for $t\in [t_0,t_0+t_\ast)$; so
$$
   \|x_t(t_0,\phi,f_\delta)\|_X 
   \leq M\|\phi\|_X+\int_{t_0}^t M \zeta_\ast(\delta)\|x_s(t_0,\phi,f_\delta)\|_X ds, 
$$
where $M:=\sup_{0\leq t\leq \tau_0}\|T(t)\|_{{\cal L}(X)}$. 
By Gronwall's inequality 
$$
   \| x_t(t_0,\phi,f_\delta) \|_X 
   \leq M\|\phi\|_X e^{M\zeta_\ast(\delta)(t-t_0)}
   \leq M\|\phi\|_X e^{M\zeta_\ast(\delta)\tau_0}, 
        \quad t\in [t_0,t_0+t_\ast),
$$
which means that $x_t(t_0,\phi,f_\delta)$ can be defined on the interval 
$[t_0,t_0+t_\ast]$ and therefore on  $[t_0,t_0+\tau_0]$ 
(cf.\,\cite[Corollary 1]{mur}). Thus it turns out that if $\|\phi\|_X$ 
is small enough, $x_t(t_0,\phi,f_\delta) $ exists on $[t_0,t_0+\tau_0]$ 
and moreover belongs to $\Omega$ in this interval. The claim readily 
follows from the fact that 
$\inf\{ W(\phi): \phi\in \Omega, \|\phi\|_X\geq r\}>0$ for small $r>0$, 
together with the property of $\Omega$.

\vskip 2.5mm

Now given $t_0\in {\bb R}^+$ and $\phi\in X$ with $W(\phi)\leq c_0$, 
consider the solution $x(t):=  x(t;t_0,\phi,f_\delta)$. 
By Proposition \ref{Proposition D2} (i) 
$$
   \| \Lambda _{\ast,\delta}(\Pi^c x_t)(s)\|_X 
   \leq \| \Lambda _{\ast,\delta}(\Pi^c x_t)\|_{Y_\eta}e^{\eta|s|} 
   \leq e^{\eta|s|}2C\|\Pi^c x_t\|_X, \quad s\in {\bb R};
$$
hence taking account of 
$\Lambda _{\ast,\delta}(\Pi^c x_t)(s)=x_{t+s}(t,\hat x_t,f_\delta)$ 
for $s\in {\bb R}$ (Proposition \ref{Proposition D3} (ii)), we get 
$$
   \|x_{t+s}(t,\hat x_t,f_\delta)\|_X 
   \leq e^{\eta \tau_0}2C\|\Pi^c x_t\|_X,  \quad s\in [-\tau_0,0],
$$ 
where $\hat{x}_t:=\Pi^cx_t+F_{\ast,\delta}(\Pi^cx_t)$. 
Set $y^\circ(t+s;t):=\lp \Psi_c, x_{t+s}(t,\hat x_t,f_\delta) \rp$. Then 
\begin{align*}
   |y^\circ(t+s;t)| 
   &\leq K_\infty \|\Psi_c\|\| x_{t+s}(t,\hat x_t,f_\delta) \|_X \\
   &\leq K_\infty \|\Psi_c\| e^{\eta \tau_0}2C\|\Pi^c x_t\|_X \\
   &\leq K_\infty \|\Psi_c\| e^{\eta \tau_0}2C a_0 \\
   &= a/2, \quad s\in [-\tau_0,0], 
\end{align*} 
hence $ y^\circ(s;t)\in S_{a/2}$ and thus $V(y^\circ(s;t))$ is well-defined for $s\in [t_0,t]$ with $t\in [t_0,t_0+\tau_0]$.  

We next confirm:  

\begin{Claim}
$\sup\{W(x_t): t\in [t_0,t_0+\tau_0]\}\leq Q$ and 
$W(x_{t_0+\tau_0}(t_0,\phi,f_\delta))\leq c_0/2$. 
\end{Claim}

Indeed, fix a $t\in [t_0,t_0+\tau_0]$ and set $z(s):=y^\circ(s;t)$ 
for $s\in [t_0,t]$. Since 
$$
   y^\circ(s;t)
   = \lp \Psi_c, x_{s}(t,\hat x_t,f_\delta) \rp
   = \lp \Psi_c, \Pi^cx_{s}(t,\hat x_t,f_\delta) \rp, \quad s\in [t_0,t],
$$
$z(s)$ is a solution of ($CE_\delta$) on $[t_0,t]$ satisfying 
$z(t)=y^\circ(t;t)=\lp\Psi_c, \Pi^c x_t\rp$. By the property (i), 
$\dot V(z(s))\leq -V(z(s))$ for $s\in [t_0,t]$, which means  
$$
   \frac{d}{ds}\big(e^{s-t}V(z(s))\big)
   = e^{s-t}\big(V(z(s))+\dot V(z(s))\big)\leq 0,
$$
so that 
$$
   V(z(t))-e^{t_0-t}V(z(t_0))
   \leq \int_{t_0}^t  \frac{d}{ds}\big(e^{s-t}V(z(s))\big)ds
   \leq 0;
$$
consequently,
\begin{align*}
   V(\lp\Psi_c, \Pi^c x_t\rp) 
   &\leq e^{t_0-t}V(y^{\circ}(t_0;t)) \\
   &= e^{t_0-t}V(\lp\Psi_c, \Pi^c x_{t_0}\rp)+ e^{t_0-t}\big( V(y^{\circ}(t_0;t))- V(\lp\Psi_c, \Pi^c x_{t_0}\rp)\big) \\
   &\leq e^{t_0-t}V(\lp\Psi_c, \Pi^c x_{t_0}\rp)+ e^{t_0-t} \bar K \big| y^{\circ}(t_0;t)- \lp\Psi_c, \Pi^c x_{t_0}\rp \big| \\
   &= e^{t_0-t}V(\lp\Psi_c, \Pi^c \phi\rp)+ e^{t_0-t} \bar K 
      \big|\lp \Psi_c, \Pi^c x_{t_0}(t,\hat x_t,f_\delta) - \Pi^c x_{t_0}\rp\big| \\
   &\leq e^{t_0-t}V(\lp\Psi_c, \Pi^c \phi\rp)+ e^{t_0-t} \bar K K_\infty\|\Psi_c\|  
                \|\Pi^c x_{t_0}(t,\hat x_t,f_\delta) - \Pi^c x_{t_0}\|_X  \\
   &\leq e^{t_0-t}V(\lp\Psi_c, \Pi^c \phi \rp)+ e^{t_0-t} \bar K K_\infty \|\Psi_c\|
           K\int_{t_0}^t e^{\mu'(\theta-t_0)}\|\xi(\theta)\|_Xd\theta,      
\end{align*}
where the last inequality is due to Proposition \ref{Proposition D4} (ii). 
Therefore, applying Proposition \ref{Proposition D6}, 
\begin{align}\label{450}
   W(x_t) 
   &= V(\lp\Psi_c, \Pi^c x_t\rp) +P\|\xi(t)\|_X \notag \\
   &\leq e^{t_0-t}V(\lp\Psi_c, \Pi^c \phi\rp) + e^{t_0-t} \bar K K_\infty \|\Psi_c\|
           K\int_{t_0}^t e^{\mu'(\theta-t_0)}\big( C\|\xi(t_0)\|_Xe^{-\beta_0(\theta-t_0)} \big) d\theta  \notag \\
   &\quad +PC\|\xi(t_0)\|_Xe^{-\beta_0(t-t_0)} \notag \\
   &\leq e^{t_0-t}V(\lp\Psi_c, \Pi^c \phi\rp) 
          + \frac{\bar K K_\infty KC\|\Psi_c\|}{\beta_0-\mu'}\|\xi(t_0)\|_X e^{t_0-t}  
          + PC\|\xi(t_0)\|_Xe^{-\beta_0(t-t_0)}. \qquad 
\end{align}
In particular, 
\begin{align*}
   W(x_{t_0+\tau_0}) 
   &\leq e^{-\tau_0}V(\lp\Psi_c, \Pi^c \phi\rp) 
         + \frac{\bar K K_\infty KC\|\Psi_c\|}{\beta_0-\mu'}\|\xi(t_0)\|_X 
         e^{-\tau_0}  
         + PC\|\xi(t_0)\|_Xe^{-\beta_0 \tau_0}  \\
   &\leq (1/2) V(\lp\Psi_c, \Pi^c \phi\rp)+ (1/4)P\|\xi(t_0)\|_X+ (1/4)P\|\xi(t_0)\|_X \\
   &= (1/2)W(x_{t_0}) 
    = (1/2)W(\phi) 
    \leq (1/2)c_0. 
\end{align*}
Since $\|\xi(t_0)\|_X\leq \|\Pi^s\phi\|_X+\|F_{\ast,\delta}(\Pi^c\phi)\|_X$, 
(\ref{450}) implies also
\begin{align*}
   &\sup\{ W(x_t): t\in [t_0,t_0+\tau_0]\} \\
   &\leq V(\lp\Psi_c, \Pi^c \phi\rp)  
         + \frac{\bar K K_\infty KC\|\Psi_c\|}{\beta_0-\mu'}\|\xi(t_0)\|_X + PC\|\xi(t_0)\|_X  \\
   &\leq V(\lp\Psi_c, \phi\rp)+\left( PC+\frac{\bar KK_\infty\|\Psi_c\|KC}{\beta_0-\mu'}\right)\big(\|\Pi^s\phi\|_X+\|F_{\ast,\delta}(\Pi^c\phi)\|_X\big) 
   = Q,   
\end{align*}
as required.

\vskip 1mm
By Claim 2, combined with Claim 1, $x(t)=x(t;t_0,\phi,f_\delta)$ is defined on $[t_0,t_0+2\tau_0]$, and 
$y^\circ(s;t)\in S_{a/2}$ still holds for $s\in [t_0,t]$ with $t\in [t_0,t_0+2\tau_0]$. 

Then we have:

\begin{Claim}
$\sup\{W(x_t): t\in [t_0+\tau_0,t_0+2\tau_0]\} \leq Q/2$ 
and  $W(x_{t_0+2\tau_0})\leq c_0/2^2$. 
\end{Claim}
 
Indeed, let $t\in [t_0+\tau_0,t_0+2\tau_0]$. By the same reasoning as 
in Claim 2 the inequality  
$$
   \frac{d}{ds}\big( e^{s-t}V(y(s;t)) \big)\leq 0, \quad s\in [t-\tau_0,t]
$$
holds; so that $V(z(t))-e^{-\tau_0}V(z(t-\tau_0))\leq 0$ and hence 
\begin{align*}
   V(\lp \Psi_c,\Pi^c x_t\rp)
   &\leq e^{-\tau_0}V(\lp\Psi_c, \Pi^c x_{t-\tau_0} \rp)
         +e^{-\tau_0} \bar K \big|\lp \Psi_c, \Pi^c 
          x_{t-\tau_0}(t,\hat x_t,f_\delta) - \Pi^c x_{t-\tau_0}\rp\big| \\
   &\leq e^{-\tau_0}V(\lp\Psi_c, \Pi^c x_{t-\tau_0}\rp)
         +e^{-\tau_0} \bar K K_\infty\|\Psi_c\|  
         \|\Pi^c x_{t-\tau_0}(t,\hat x_t,f_\delta) - \Pi^c x_{t-\tau_0}\|_X \\
   &\leq e^{-\tau_0}V(\lp\Psi_c, \Pi^c x_{t-\tau_0} \rp)
         +e^{-\tau_0} \bar K K_\infty \|\Psi_c\|
         K\int_{t-\tau_0}^t e^{\mu'(\theta-t+\tau_0)}\|\xi(\theta)\|_Xd\theta. 
\end{align*}
Therefore, correspondingly to (\ref{450}),
\begin{align*}
   W(x_t) 
   &\leq e^{-\tau_0}V(\lp\Psi_c, \Pi^c x_{t-\tau_0} \rp)   \\
   &\quad +e^{-\tau_0} \bar K K_\infty \|\Psi_c\|
          K\int_{t-\tau_0}^t e^{\mu'(\theta-t+\tau_0)} \big(C\|\xi(t-\tau_0)\|_X          e^{-\beta_0(\theta-t+\tau_0)}\big) d\theta \\
   &\quad +PC\|\xi(t-\tau_0)\|_Xe^{-\beta_0 \tau_0}   \\
   &\leq e^{-\tau_0}V(\lp\Psi_c, \Pi^c x_{t-\tau_0} \rp) 
         +\frac{\bar K K_\infty KC\|\Psi_c\|}{\beta_0-\mu'}
         \|\xi(t-\tau_0)\|_X e^{-\tau_0} \\ 
  &\quad +PC\|\xi(t-\tau_0)\|_Xe^{-\beta_0 \tau_0}  \\
  &\leq (1/2) V(\lp\Psi_c, \Pi^c x_{t-\tau_0} \rp)
        + (1/4)P\|\xi(t-\tau_0)\|_X+ (1/4)P\|\xi(t-\tau_0)\|_X \\
  &= (1/2)W(x_{t-\tau_0})  \\
  &\leq (1/2)\sup \{W(x_\tau): \tau\in [t_0,t_0+\tau_0]\}
   \leq Q/2,
\end{align*}
where the last inequality follows from Claim 2. Letting $t=t_0+2\tau_0$ 
in the above, we also see from claim 2 that
$$
   W(x_{t_0+2\tau_0})
   \leq \frac{1}{2}W(x_{t_0+\tau_0})
   \leq \frac{1}{2}\cdot \frac{c_0}{2}
   = \frac{c_0}{2^2},
$$
and Claim 3 holds.

\vskip 2.5mm

Repeating this argument, one can deduce in general that 
$x(t)=x(t;t_0,\phi,f_\delta)$ is defined on $[t_0,t_0+n\tau_0]$, and 
$y^\circ(s;t)\in S_{a/2}$ holds for $s\in [t_0,t]$ with 
$t\in [t_0,t_0+n\tau_0]$ for any $n\in {\bb N}$. Moreover, 
$$
   \sup\{W(x_t): t\in [t_0+(n-1)\tau_0,t_0+n\tau_0]\}
   \leq \frac{Q}{2^{n-1}} 
        \quad \text{and} \quad W(x_{t_0+n\tau_0})\leq \frac{c_0}{2^n}
$$ 
for $n\in {\bb N}$. This means that $x(t)=x(t;t_0,\phi,f_\delta)$ is 
actually defined on $[t_0,\infty)$ and that 
$$
   V(\lp \Psi_c,x_t(t_0,\phi,f_\delta)\rp)
   + P\|\Pi^sx_t-F_{\ast,\delta}(\Pi^c x_t)\|_X
   \leq Q\,2^{-(t-t_0)/\tau_0}, \quad t\in [t_0,\infty).
$$
In view of (i) and $P>1$, 
$$
   b(|\lp \Psi_c,x_t(t_0,\phi,f_\delta)\rp|)
   \leq Q\,2^{-(t-t_0)/\tau_0}\leq b(a), 
   \quad \|\Pi^sx_t-F_{\ast,\delta}(\Pi^c x_t)\|_X \leq Q\,2^{-(t-t_0)/\tau_0}.
$$
Since
$$
  \|\Pi^c x_t(t_0,\phi,f_\delta)\|_X 
  = \|\Phi_c \lp \Psi_c,x_t(t_0,\phi,f_\delta)\rp\|_X
  \leq \|\Phi_c\|\,b^{-1}\big(Q\,2^{-(t-t_0)/\tau_0}\big)
$$
with $\|\Phi_c\|:= \big(\sum_{j=1}^{d_c}\|\phi_j\|_X^2\big)^{1/2}$ and
\begin{align*}
   \| \Pi^s x_t(t_0,\phi,f_\delta)\|_X
   &\leq \|\Pi^s x_t-F_{\ast,\delta}(\Pi^c x_t)\|_X 
         + \|F_{\ast,\delta}(\Pi^c x_t)\|_X \\
   &\leq Q\,2^{-(t-t_0)/\tau_0}+L(\delta)\|\Pi^c x_t\|_X,
\end{align*}
so we obtain that for any $\phi \in \Omega$ and $t\in [t_0,\infty)$
\begin{align*}
   \|x_t(t_0,\phi,f_\delta)\|_X
   &\leq \|\Pi^c x_t\|_X + \|\Pi^s x_t\|_X  \\
   &\leq Q\,2^{-(t-t_0)/\tau_0}
         + (1+L(\delta))\|\Phi_c\|\,b^{-1}\big(Q\,2^{-(t-t_0)/\tau_0}\big),
\end{align*}
which shows that the zero solution of ($E_\delta$) is uniformly asymptotically stable. \qed

\vskip 5mm


Before concluding this section, we will provide an example to illustrate 
how our Theorem \ref{Theorem D3} is available for stability analysis of 
some concrete equations. Let us consider nonlinear (scalar) integral equation 
\begin{align}\label{1000}
   x(t)=\nu\int_{-\infty}^tP(t-s)x(s)ds+f(x_t),
\end{align}
where $\nu$ is a nonnegative real parameter, $P$ is a nonnegative 
continuous function on ${\bb R}^+$ satisfying $\int_0^{\infty}P(t)dt=1$ 
together with the condition 
$\|P\|_{1,\rho}:=\int_0^{\infty}P(t)e^{\rho t}dt<\infty$ and 
$\|P\|_{\infty,\rho}:=\mathrm{ess}\,\sup\{ P(t)e^{\rho t}: t\geq 0\}<\infty$ 
for some positive constant $\rho$, and $f\in C^1(X;{\bb C})$, 
$X:=L^1_{\rho}({\bb R}^-;{\bb C})$, satisfies $f(0)=0$ and $Df(0)=0$.  
Eq. (\ref{1000}) is written as Eq. ($E$) with $m=1$ and $K\equiv \nu P$. 
The characteristic operator $\Delta(\lambda)$ of Eq. (\ref{1000}) 
is given by 
$\Delta(\lambda)=1-\nu \int_0^{\infty}P(t)e^{-\lambda t}dt$. 
It is easy to see that if $\nu>1$, then $\Delta(\lambda_0)=0$ for some 
positive $\lambda_0$; hence $\Sigma^u\neq\emptyset$. Observe that 
$|\Delta(\lambda)|\geq 1-\nu$ if $\mathrm{Re}\,\lambda \geq 0$ and 
$0\leq \nu\leq 1$. Thus, if $0\leq \nu<1$, then 
$\Sigma^u\cup\Sigma^c=\emptyset$. Hence, by virtue of 
{\bf the principle of linearized stability} for integral equations 
(e.g., \cite[Theorem 3.15]{die-gyl}), we get the following:  

\begin{Prop}\label{Proposition E1}
Under the above conditions on Eq. (\ref{1000}), the following statements 
hold true;
\begin{enumerate}
   \item if $0\leq \nu<1$, then the zero solution of  Eq. (\ref{1000}) 
         is exponentially stable (in $L^1_{\rho}$);
   \item if $\nu>1$, then the zero solution of  Eq. (\ref{1000}) 
         is unstable (in $L^1_{\rho}$)
\end{enumerate}
\end{Prop}

In the remainder of this section, we will treat Eq. (\ref{1000}) in 
the critical case $\nu=1$, and investigate stability property for the 
zero solution of Eq. (\ref{1000}) by applying Theorem \ref{Theorem D3}. 
In case $\nu=1$, we easily see that $\Sigma^u=\emptyset$ and $\Sigma^c=\{0\}$. 
Indeed, in this case, $0$ is a simple root of the equation 
$\Delta(\lambda)=0$, and $E^c$ is $1$-dimensional space with 
a basis $\{\phi_1\}$, $\phi_1\equiv 1$, together with 
$\{\psi_1\},\ \psi_1\equiv 1$, as the dual basis of $\{\phi_1\}$; 
see \cite{murmat} for details.  The projection $\Pi^{c}$ is given 
by the formula $\Pi^c\phi=\Phi_c\lp \Psi_c,\phi \rp, \ \forall\phi\in X$, 
and hence 
\begin{align*}
   \Pi^c\phi 
    = \phi_1\lp\psi_1,\phi \rp
   &= \phi_1\left(\int_{-\infty}^0\int_{\theta}^0
      \psi_1(\xi-\theta)P(-\theta)\phi(\xi)d\xi d\theta\right)\\
   &= \Phi_c\left(\int_{-\infty}^0P(-\theta)
      \left(\int_{\theta}^0\phi(\xi)d\xi\right)d\theta\right).
\end{align*}
Thus, for a solution $x(t)$ of Eq. (\ref{1000}), the component $z_c(t)$ 
of $\Pi^cx_t$ with respect to $\Phi_c$ is given by
$$
   z_c(t)=\int_{-\infty}^t\hat{P}(t-s)x(s)ds
$$
with $\hat{P}(t):=\int_t^{\infty}P(\tau)d\tau$, because of 
\begin{align*}
   z_c(t)
   &= \int_{-\infty}^0P(-\theta)\left(\int_{\theta}^0x(t+\xi)d\xi\right)d\theta 
    = \int_{-\infty}^0P(-\theta)\left(\int_{t+\theta}^tx(s)ds\right)d\theta \\
   &= \int_{-\infty}^tP(t-\tau)\left(\int_{\tau}^tx(s)ds\right)d\tau 
    = \int_{-\infty}^t\left(\int_{t-s}^{\infty}P(w)dw\right)x(s)ds. 
\end{align*}
Observe that $z_c(t)$ satisfies the ordinary equation 
$$
   \dot{z}_c(t)= \hat{P}(0)x(t)+\int_{-\infty}^t(-P(t-s))x(s)ds
   =x(t)-\int_{-\infty}^tP(t-s)x(s)ds,
$$
that is, 
$
   \dot{z}_c(t)=f(x_t)=f(\Phi_cz_c(t)+\Pi^sx_t). 
$
In particular, if $x$ is a solution of Eq. (\ref{1000}) satisfying 
$x_t\in W^c_{\rm loc}(r,\delta)$ on an interval $J$, then   
$\Pi^sx_t=F_{\ast}(\Phi_cz_c(t))$ on $J$; hence we get  
$
   \dot{z}_c(t)=f(\Phi_cz_c(t)+F_{\ast}(\Phi_cz_c(t)))
$
on $J$.  
This observation leads to that $G_c=0$ and $H_c=1$ in the central equation 
($CE$); in fact, by noticing that $\Sigma^c=\{0\}$ and 
$H_cx=\lim_{n\to\infty}\lp\psi_1,\Gamma^nx \rp =x,\ \forall x\in {\bb C}$, 
one can also certify this fact.  
Let us assume that $f$ is of the form 
\begin{align}\label{1002}
   f(\phi)
   = \epsilon\bigg(\int_{-\infty}^0\hat{P}(-\theta)\phi(\theta)d\theta\bigg)^3
     +g(\phi), \quad \forall \phi\in X,
\end{align}
where $\epsilon$ is a nonzero real number, and $g\in C^1(X;{\bb C})$ 
satisfies $|g(\phi)|=o(\|\phi\|^3_X)$ as $\|\phi\|_X\to 0$ (here, $o$ 
means Landau's notation "small oh"). Indeed, noting that 
$\|\hat{P}\|_{\infty,\rho}\leq (1/\rho)\|P\|_{\infty,\rho}$ and 
$\|\hat{P}\|_{1,\rho}\leq (1/\rho)\|P\|_{1,\rho}$, we see that the 
function $f$ given by (\ref{1002}) satisfies $f\in C^1(X;{\bb C})$ 
and $f(0)=Df(0)=0$. Since 
\begin{align*}
   f(\Phi_cz_c(t)+F_{\ast}(\Phi_cz_c(t)))
   &= f(x_t)=\epsilon\left(\int_{-\infty}^0\hat{P}(-\theta)x_t(\theta)d\theta\right)^3+g(x_t) \\
   &= \epsilon\left(\int_{-\infty}^t\hat{P}(t-s)x(s)ds\right)^3+g(x_t)\\
   &= \epsilon(z_c(t))^3+g(\Phi_cz_c(t)+F_{\ast}(\Phi_cz_c(t))),
\end{align*}
we get $f(\Phi_cw+F_{\ast}(\Phi_cw))=\epsilon w^3+
g(\Phi_c w+F_{\ast}(\Phi_c w)), \forall w\in {\bb C}.$ 
Thus, the central equation ($CE$) of Eq. (\ref{1000}) becomes to the equation 
\begin{align}\label{1001}
   \dot{z}=\epsilon z^3+g(\phi_1 z+F_{\ast}(\phi_1z)). 
\end{align}
Since $F_{\ast}(\phi_1z)=o(z)$ as $z\to 0$ by Theorem \ref{Theorem D1a}-(i), 
it follows that $\dot{z}=\epsilon z^3+o(z^3)$ as $z\to 0$; consequently, 
one can easily see that the zero solution of Eq. (\ref{1001}) is uniformly 
asymptotically stable (resp. unstable) if $\epsilon<0$ (resp. $\epsilon>0$).
Therefore, by virtue of Theorem \ref{Theorem D3}, we get the following result:

\begin{Prop}\label{Proposition D8} 
Let $\nu=1$ in Eq. (\ref{1000}), and assume that $f$ is of 
the form (\ref{1002}) (with nonzero constant $\epsilon$ and 
$g(\phi)=o(\|\phi\|^3_X)$ as $\|\phi\|_X\to 0$ with $g\in C^1(X;{\bb C})$). 
Then 
\begin{enumerate}
   \item if $\epsilon<0$, then the zero solution of Eq. (\ref{1000}) 
         is uniformly asymptotically stable (in $L^1_{\rho}$);
   \item if $\epsilon>0$, then the zero solution of Eq. (\ref{1000}) 
         is unstable (in $L^1_{\rho}$).
\end{enumerate} 
\end{Prop}

\medskip


\section{Appendix}

In this appendix we will prove the $C^1$-smoothness of the center 
manifold $W^c_\delta$ of the equilibrium point $0$ of Eq.($E_\delta$), 
and give stable/unstable manifold theorems (of the equilibrium point $0$) 
of Eq.($E$).  

\subsection{Smoothness of the center manifold $W^c_\delta$}

For a Banach space $U$ with norm $\|\cdot\|_{U}$ and $\lambda\geq 0$, 
we define    
$$
   BC^\lambda({\bb R}; U)
   = \big\{y\in C({\bb R};U): 
     \sup_{t\in {\bb R}}\|y(t)\|_{U}\,e^{-\lambda |t|}<\infty \big\}.
$$
$BC^{\lambda}({\bb R}; U)$ is a Banach space normed with 
$\|y\|_{BC^{\lambda}({\bb R}; U)}:=
\sup_{t\in {\bb R}}\|y(t)\|_{U}\,e^{-\lambda |t|}$. 
We use, for abbreviation, the notation $\|\cdot\|_{\lambda,\,U}$ 
instead of $\|\cdot\|_{BC^{\lambda}({\bb R}; U)}$. 
Evidently, if  $\lambda\leq \lambda'$, there is an inclusion map 
$BC^{\lambda}({\bb R}; U)\hookrightarrow BC^{\lambda'}({\bb R}; U)$ with 
$$
   \|y\|_{\lambda',\,U}\leq \|y\|_{\lambda,\,U}  
   \quad \text{for} \enskip y\in  BC^{\lambda}({\bb R}; U).
$$
In what follows, for $\nu\geq 0$ we denote the inclusion map from 
$BC^{\lambda}({\bb R};U)$ to $BC^{\lambda+\nu}({\bb R};U)$ by the 
same notation $J_{\nu}$ for all $\lambda\geq 0$ and any Banach space $U$. 
Clearly, $J_{\nu}$ belongs to 
${\cal L}(BC^{\lambda}({\bb R};U); BC^{\lambda+\nu}({\bb R};U))$. 

By restricting the functions with values in an open subset ${\cal O}$ 
of $U$, we also  use the symbols $BC^{\lambda}({\bb R};{\cal O})$, 
$BC^{\lambda'}({\bb R};{\cal O})$ and so on to denote open subsets 
of the spaces $BC^{\lambda}({\bb R};U)$, $BC^{\lambda'}({\bb R};U)$ 
and so on, respectively. 

Let ${\cal O}$ be an open set of $U$, and consider a continuous map 
$h:{\cal O}\to V$. Then $h$ induces a map $\widetilde{h}$ from  
$C({\bb R}; {\cal O})$ to $C({\bb R}; V)$ by letting 
$$ 
   \big[\widetilde{h}(y)\big](t):=h(y(t)) 
   \quad \text{for} \enskip y\in C({\bb R};{\cal O})
   \enskip \text{and} \enskip t\in {\bb R}. 
$$
We recall that if $\sup_{u\in {\cal O}}\|h(u)\|_V < \infty$, 
the induced map $\widetilde h$ is continuous as a map from 
$BC^{\lambda }({\bb R};{\cal O})$ to $BC^{\lambda'}({\bb R};V)$ 
for $\lambda'>0$ (cf.\,\cite[Appendix IV]{diek}).

\begin{Prop}\label{Proposition D10}
$\Lambda_{\ast,\delta}$ is of class $C^1$ as a map from 
$E^c$ to $Y_{\eta'}$ for $\eta'\in (\eta,\alpha)$. 
\end{Prop} 

\noindent
{\it Proof.} Let $\hat\eta$ be a positive number such that 
$\eta<\hat\eta<\alpha$ and set
$$
   C_\ast(\lambda )
   := \zeta_\ast(\delta_1)CC_1
      \left(\frac{1}{\lambda-\varepsilon }+\frac{2}{\alpha +\lambda }
      +\frac{2}{\alpha -\lambda }\right) 
      \quad \text{for} \enskip \lambda \in [\eta, \hat\eta].
$$
By taking $\delta_1>0$ small if necessary, we may assume that 
$$
   C_\ast(\lambda) <\frac{1}{2}, \quad \zeta_\ast(\delta)<1 
   \quad \text{for} \enskip 0<\delta\leq\delta_1 
   \enskip \text{\and} \enskip \lambda \in [\eta, \hat\eta].
$$
Now let $\lambda\in [\eta,\hat\eta]$, $\mu$ be a nonnegative number 
with $\lambda +\mu\leq \hat\eta$ and  
$v\in BC^\mu({\bb R}; {\cal L}(X;{\bb C}^m))$. 
Consider a map  
${\cal H}_{0}(v) : Y_{\lambda }\to Y_{\lambda +\mu}$ defined by 
\begin{align*}
   \big[{\cal H}_{0}(v)w_{0}\big](t)
   &:= \lim_{n\to\infty}\int_{0}^{t}T^c(t-s)\Pi^{c}\Gamma^{n} v(s)w_{0}(s)ds \\
   &\quad -\lim_{n\to\infty}\int_{t}^{\infty}T^u(t-s)\Pi^{u}\Gamma^{n} v(s)w_{0}(s)ds \\ 
   &\quad +\lim_{n\to\infty}\int_{-\infty}^{t}T^s(t-s)\Pi^{s}\Gamma^{n} v(s)w_{0}(s)ds, \quad t\in {\bb R}
\end{align*}
for $w_{0}\in Y_{\lambda }$.  By (\ref{286a5}) it follows that 
\begin{align*}
   \| [{\cal H}_0(v)w_0](t)\|_X 
   &\leq \Big| \int_{0}^{t}CC_1\zeta_{\ast}(\delta)e^{\varepsilon |t-s|}
         \big(\|v\|_{\mu,\, {\cal L}(X;{\bb C}^m)}e^{\mu|s|}\big)
         \big(\|w_0\|_{Y_{\lambda }}e^{\eta |s|}\big) ds \Big| \\
   &\quad +\int_{t}^{\infty}CC_1\zeta_{\ast}(\delta)e^{\alpha (t-s)}
         \big(\|v\|_{\mu,\, {\cal L}(X;{\bb C}^m)}e^{\mu|s|}\big)
         \big(\|w_0\|_{Y_{\lambda }}e^{\eta |s|}\big) ds \\
   &\quad +\int_{-\infty}^{t}CC_1 \zeta_{\ast}(\delta)e^{-\alpha(t-s)}
         \big(\|v\|_{\mu,\, {\cal L}(X;{\bb C}^m)}e^{\mu|s|}\big)
         \big(\|w_0\|_{Y_{\lambda }}e^{\eta|s|}\big) ds \\
   &\leq C_\ast(\lambda +\mu)\|v\|_{\mu,\, {\cal L}(X;{\bb C}^m)} 
         \|w_{0}\|_{Y_{\lambda }}e^{(\lambda +\mu)|t|}
\end{align*}
for any $w_0\in Y_{\lambda}$ and $t\in {\bb R}$; hence 
\begin{align}\label{525}
   \|{\cal H}_{0}(v)w_{0}\|_{Y_{\lambda +\mu}}
   \leq \frac{1}{2} \|v\|_{\mu,\, {\cal L}(X;{\bb C}^m)} 
        \|w_{0}\|_{Y_{\lambda }}, \quad w_{0}\in Y_{\lambda }.  
\end{align}
${\cal H}_{0}(v)$ induces a bounded linear map from 
$Y_{\lambda }^{(1)}={\cal L}(E^c;Y_\lambda )$ to 
$Y_{\lambda +\mu}^{(1)}={\cal L}(E^c;Y_{\lambda +\mu})$, denoted 
${\cal H}(v)$, by 
$$
   \big[[{\cal H}(v)w](t)\big]\phi
   := [{\cal H}_{0}(v)w\phi](t) 
      \quad \text{for} \enskip w\in Y_{\lambda }^{(1)}
      \enskip \text{and} \enskip \phi \in E^{c}, 
$$
where $w\phi\in Y_{\lambda }$ is given by $[\,w\phi\,](t):=w(t)\phi$ 
for $t\in {\bb R}$, and (\ref{525}) yields the  estimate    
\begin{align}\label{526}
   \|{\cal H}(v)w\|_{Y_{\lambda +\mu}^{(1)}}
   \leq \frac{1}{2} \|v\|_{\mu,\, {\cal L}(X;{\bb C}^m)} 
        \|w\|_{{Y_\lambda }^{(1)}}, \quad w\in Y_{\lambda }^{(1)}.  
\end{align}
By setting  $\lambda =\eta$, $\mu=0$ and 
$v_\ast(\psi):=\widetilde{Df_{\delta}}(\Lambda_{\ast,\delta}(\psi))$ for 
$\psi\in E^c$, we can consider a linear equation in $Y_{\eta}^{(1)}$
\begin{align}\label{552}
   A_{1}=T^{c}(\cdot)+{\cal H}(v_\ast(\psi) ) A_{1}, 
   \quad \psi \in E^c
\end{align}
since $T^{c}(\cdot)$ belongs to $Y_{\eta}^{(1)}$. Because 
$\| \widetilde{Df_{\delta}}(\Lambda_{\ast,\delta}(\psi)) \|_{0,\,{\cal L}(X;{\bb C}^m)}\leq \zeta_{\ast}(\delta )$ (cf.\,(\ref{351})),  
(\ref{526}) implies 
\begin{align}\label{552a}
   \| {\cal H}(v_\ast(\psi)) \|_{{\cal L}(Y_{\eta}^{(1)})} 
   \leq \frac{1}{2} \zeta_{\ast}(\delta )
   \leq \frac{1}{2},
\end{align}
and hence   
$A_{1}=A_{1}(\psi)$ is uniquely determined for each $\psi \in E^{c}$ 
and is given by
\begin{align}\label{553}
   A_{1}(\psi)
   = \big( I_{Y_{\eta}^{(1)}} - {\cal H}(v_\ast(\psi) )\big)^{-1} T^{c}(\cdot)
   = \sum_{n=0}^{\infty} ({\cal H}( v_\ast(\psi) ))^n T^{c}(\cdot).
\end{align}

Let us take any $\eta' \in (\eta,\hat{\eta}]$.  We will verify that 
$J_{\eta'-\eta}\Lambda_{\ast,\delta }(\psi)$ is differentiable and 
\begin{align}\label{553b}
   D( J_{\eta'-\eta}\Lambda_{\ast,\delta }(\psi) )
   = J_{\eta'-\eta}A_1(\psi), \quad \psi\in E^c
\end{align}
holds. Notice that ${\cal F}_{\delta}(\psi+h, \Lambda_{\ast, \delta}(\psi+h))-
{\cal F}_{\delta}(\psi, \Lambda_{\ast, \delta}(\psi+h))=T^c(\cdot)h=A_1(\psi)h-{\cal H}_0(v_{\ast}(\psi))A_1(\psi)h$, $\forall h\in E^c$, 
by (\ref{287}) and (\ref{552}). We thus get  
\begin{align*} 
   &\Lambda_{\ast,\delta}(\psi+h)-\Lambda_{\ast,\delta}(\psi)-A_{1}(\psi)h \\
   &= {\cal F}_{\delta}(\psi+h, \Lambda_{\ast,\delta}(\psi+h))
      -{\cal F}_{\delta}(\psi, \Lambda_{\ast,\delta}(\psi))- A_1(\psi)h \\ 
   &= {\cal F}_{\delta}(\psi, \Lambda_{\ast,\delta}(\psi+h))
      -{\cal F}_{\delta}(\psi, \Lambda_{\ast,\delta}(\psi))
      -{\cal H}_{0}(v_\ast(\psi) )A_{1}(\psi)h
\end{align*}
for $h\in E^{c}$. Since 
$$
   f_{\delta}(\Lambda_{\ast,\delta}(\psi+h)(s))-f_{\delta}(\Lambda_{\ast,\delta}(\psi)(s)) 
   = \int_{0}^{1}Df_{\delta}(\Lambda_{\sigma}(\psi)(s))d\sigma\,\big( \Lambda_{\ast,\delta}(\psi+h)(s)-\Lambda_{\ast,\delta}(\psi)(s)\big),
$$
where $\Lambda_{\sigma}(\psi):=(1-\sigma)\Lambda_{\ast,\delta}(\psi)+\sigma \Lambda_{\ast,\delta}(\psi+h)$ for $\sigma\in [0,1]$, 
it follows that
$$
   {\cal F}_{\delta}(\psi, \Lambda_{\ast,\delta}(\psi+h))
      -{\cal F}_{\delta}(\psi, \Lambda_{\ast,\delta}(\psi))
   = {\cal H}_{0}(v_h(\psi)) (\Lambda_{\ast,\delta}(\psi+h)
     -\Lambda_{\ast,\delta}(\psi)),
$$
where $v_{h}(\psi)$ is an element in $BC({\bb R}; {\cal L}(X;{\bb C}^m))$ 
defined by
$$
   [v_h(\psi)](s):
   =\int_{0}^{1}Df_{\delta}( \Lambda_{\sigma}(\psi)(s)) d\sigma
    \quad \text{for} \enskip s\in {\bb R}.
$$ 
Hence 
\begin{align*}
   \Lambda_{\ast,\delta}(\psi+h)-\Lambda_{\ast,\delta}(\psi)-A_{1}(\psi)h 
   &= {\cal H}_{0} ( v_h(\psi) -v_\ast(\psi))
      (\Lambda_{\ast,\delta}(\psi+h)-\Lambda_{\ast,\delta}(\psi)) \\
   &\quad +{\cal H}_{0}(v_\ast(\psi) )(\Lambda_{\ast,\delta}(\psi+h)-\Lambda_{\ast,\delta}(\psi)-A_{1}(\psi)h),
\end{align*}
and, applying  (\ref{525}) and Proposition \ref{Proposition D2} (i), we get
\begin{align*}
   &\|J_{\eta'-\eta}\Lambda_{\ast,\delta}(\psi+h)
      -J_{\eta'-\eta}\Lambda_{\ast,\delta}(\psi)
      -J_{\eta'-\eta}A_{1}(\psi)h\|_{Y_{\eta'}} \\
   &\leq (1/2)\|v_h(\psi)-v_\ast(\psi) \|_{\eta'-\eta,\,{\cal L}(X;{\bb C}^m)}\,
    \|\Lambda_{\ast,\delta}(\psi+h)-\Lambda_{\ast,\delta}(\psi)\|_{Y_{\eta}} \\
   &\quad +(1/2) \| v_\ast(\psi) \|_{0,\,{\cal L}(X;{\bb C}^m)}\, 
        \| J_{\eta'-\eta}\Lambda_{\ast,\delta}(\psi+h)
          -J_{\eta'-\eta}\Lambda_{\ast,\delta}(\psi)
          -J_{\eta'-\eta}A_{1}(\psi)h \|_{Y_{\eta'}}  \\
   &\leq C \|v_h(\psi)-v_\ast(\psi)\|_{\eta'-\eta,\,{\cal L}(X;{\bb C}^m)}\,\|h\|_X  \\
   &\quad +(1/2)\, \| J_{\eta'-\eta}\Lambda_{\ast,\delta}(\psi+h)
          -J_{\eta'-\eta}\Lambda_{\ast,\delta}(\psi)
          -J_{\eta'-\eta}A_{1}(\psi)h \|_{Y_{\eta'}},
\end{align*}
so that 
\begin{align}\label{553a}
   &\|J_{\eta'-\eta}\Lambda_{\ast,\delta}(\psi+h)
      -J_{\eta'-\eta}\Lambda_{\ast,\delta}(\psi)
      -J_{\eta'-\eta}A_{1}(\psi)h\|_{Y_{\eta'}} \notag \\
   &\leq 2C \|v_h(\psi)-v_\ast(\psi)\|_{\eta'-\eta,\,{\cal L}(X;{\bb C}^m)} \, \|h\|_X.    
\end{align}
Notice that the continuous map 
$Df_{\delta}:S_{\delta}\to {\cal L}(X;{\bb C}^m)$ satisfies 
$\sup_{\phi\in S_{\delta}}\|Df_{\delta}(\phi)\|_{ {\cal L}(X;{\bb C}^m)}
\leq \zeta_{\ast}(\delta)$, and hence by the fact stated in the preceding 
paragraph of the proposition, the induced map  
$\widetilde{Df_{\delta}}:BC^{\eta }({\bb R}; S_{\delta})\to 
BC^{\eta'-\eta}({\bb R}; {\cal L}(X;{\bb C}^m))$ 
is continuous because of $\eta'-\eta>0$.   
Since 
$
  v_h(\psi) - v_\ast(\psi) 
  = \int_{0}^{1} \big( \widetilde{Df_{\delta}}(\Lambda_{\sigma}(\psi) )
    - \widetilde{Df_{\delta}}(\Lambda_{\ast}(\psi) ) \big)d\sigma
$
in $ BC^{\eta'-\eta}({\bb R};{\cal L}(X;{\bb C}^m))$, the following 
inequality holds true:
$$
   \| v_h(\psi) -v_\ast(\psi) \|_{\eta'-\eta,\,{\cal L}(X;{\bb C}^m)}
   \leq \sup_{0\leq \sigma\leq 1} 
        \|\widetilde{Df_{\delta}}(\Lambda_{\sigma}(\psi))
        -\widetilde{Df_{\delta}}(\Lambda_{\ast,\delta}(\psi))\|_{\eta'-\eta,\,{\cal L}(X;{\bb C}^m)} .
$$
Hence $ \lim_{h\to 0} \|v_h(\psi) -v_\ast(\psi) \|_{\eta'-\eta,\,{\cal L}(X;{\bb C}^m)}=0$
because of 
$$
  \sup_{0\leq \sigma\leq 1}\|\Lambda_{\sigma}(\psi)
       -\Lambda_{\ast,\delta}(\psi)\|_{Y_{\eta}} 
  \leq \|\Lambda_{\ast,\delta}(\psi+h)-\Lambda_{\ast,\delta}(\psi)\|_{Y_{\eta}} 
  \leq 2C\|h\|_X \to 0
$$
as $\|h\|_X\to 0$. Thus (\ref{553a}) implies the differentiability 
of the map $J_{\eta'-\eta}\Lambda_{\ast,\delta }(\psi)$ as well as 
$D(J_{\eta'-\eta}\Lambda_{\ast,\delta }(\psi)) = J_{\eta'-\eta}A_{1}(\psi)$.

We will finally certify the $C^1$-smoothness of the map 
$J_{\eta'-\eta}\Lambda_{\ast,\delta }(\psi)$. By (\ref{553}) 
$$
   J_{\eta'-\eta}A_1(\psi)
   = \sum_{n=0}^{\infty} J_{\eta'-\eta} ({\cal H}( v_\ast(\psi)) )^n T^{c}(\cdot).
$$
This series in $Y_{\eta'}^{(1)}$ converges uniformly for $\psi\in E^c$ 
because the norm of each term $J_{\eta'-\eta}({\cal H}(v_\ast(\psi)))^n$ 
does not exceed $(1/2)^n$ by (\ref{552a}).  So,  by virtue of (\ref{553b}), 
it suffices to prove the continuity of the term 
$J_{\eta'-\eta} ({\cal H}(v_\ast(\psi)))^n$ as a map from $E^c$ to 
${\cal L}(Y_{\eta}^{(1)}, Y_{\eta'}^{(1)})$. Given a positive integer $n$, 
put $a_{1}:=(\eta'-\eta)/n$. It is then easy to see that as a map from 
$Y_{\eta}^{(1)}$ to $Y_{\eta'}^{(1)}$
$$
  J_{\eta'-\eta} ({\cal H}(v_\ast(\psi)))^n 
  = J_{a_{1}}^{n} ({\cal H}( v_\ast(\psi)))^n
  = (J_{a_{1}} {\cal H} ( v_\ast(\psi)))^n    
  = ( {\cal H} ( J_{a_{1}} v_\ast(\psi)))^n
$$
holds for $\psi\in E^c$. Since 
$\sup_{\phi\in S_{\delta}}\|Df_{\delta}(\phi)\|_{{\cal L}(X;{\bb C}^m)}
\leq \zeta_{\ast}(\delta)$ again, the map 
$J_{a_{1}}\widetilde{Df_{\delta}}:BC^{\eta}({\bb R}; S_{\delta})\to 
BC^{a_1}({\bb R}; {\cal L}(X;{\bb C}^m))$ 
is continuous by the fact stated in the preceding paragraph of 
the proposition. Then $J_{a_{1}} v_\ast(\psi)$ is continuous in $\psi$; 
and so is $J_{\eta'-\eta} ({\cal H}(v_\ast(\psi)))^n$. 
This completes the proof.  \qed

\vskip 1.5mm
\begin{Rem} {\rm Let $k$ be a positive integer. Then under the assumption 
that $f$ is of class $C^k$, we can establish the $C^k$-smoothness of 
the center manifold $W^c_\delta$. In fact, the continuity of each term of 
the formal series, given by ($m$ times) term-wise differentiation of 
(\ref{553}) ($m=1,2,\ldots,k-1$), is guaranteed, together with the 
uniform convergence of the series, if we regard $\Lambda _{\ast,\delta}$ 
as a map from $E^c$ to $Y_{\hat\eta}$ for a suitable 
$\hat\eta\in (\eta,\alpha)$; for details  see \cite[Proposition 4]{murnag}. }
\end{Rem}

\subsection{Stable manifold theorem and unstable manifold theorem}  

In this subsection we will give (local) stable/unstable manifold 
theorems for the integral equation
\begin{align}
 x(t)=\int_{-\infty}^t K(t-s)x(s)ds+ f(x_{t}) \tag{$E$}
\end{align}
under the assumption that the zero solution of ($E$) is hyperbolic.

For $r>0$ and $\delta>0$ we set 
$$
   W^{s}_{{\rm loc}}(r,\delta)=\big\{\phi \in X: \| \Pi^{s}\phi \|_X < r, 
   \enskip \| x_{t}(0,\phi,f)\|_X <\delta, \ t\in {\bb R}^{+}\big\}.
$$
Then $W^{s}_{{\rm loc}}(r,\delta)$ is called the local stable manifold 
(of the equilibrium point $0$) of Eq.($E$). 

\begin{Thm}\label{Theorem D6}
Assume that the zero solution of Eq.($E$) is hyperbolic and that 
$f\in C^{k}(X;{\bb C}^m)$ with $f(0)=Df(0)=0$. Then there exist positive 
numbers $r$, $\delta$, and a $C^k$-map $F^s:B_{E^{s}}(r)\to E^{u}$ 
with $F^s(0)=0$, together with an open neighborhood $\Omega_0$ of 
$0$ in $X$, such that the following properties hold: 
\begin{enumerate}
   \item $W^{s}_{{\rm loc}}(r,\delta)={\rm graph}\,F^s$; moreover,   $W^{s}_{{\rm loc}}(r,\delta)$ is tangent 
         to $E^{s}$ at zero.
   \item To any $\beta\in (0,\alpha)$ there corresponds a positive 
         constant $M$ such that 
$$
   \| x_{t}(0,\phi,f) \|_X
   \leq Me^{-\beta t}\|\phi \|_X, 
   \quad t\in {\bb R}^{+}, \enskip \phi\in W^{s}_{{\rm loc}}(r,\delta).
$$
   \item $W^{s}_{{\rm loc}}(r,\delta)$ is locally positively invariant 
         for Eq.($E$), that is, if $\phi\in W^{s}_{{\rm loc}}(r,\delta)$, 
         we have $x_{\tau}(0,\phi,f)\in W^{s}_{{\rm loc}}(r,\delta)$ 
         for $\tau\in {\bb R}^{+}$ whenever 
         $\Pi^{s}x_{\tau}(0,\phi,f)\in B_{_{E^{s}}}(r)$.
   \item There exists a positive constant $\beta_1$ with the property 
         that if $x(t)$ is a solution of Eq.($E$) on an interval 
         $J=[t_0,t_1]$ satisfying $x_t\in \Omega_0$ on $J$, then the 
         inequality 
$$
   \|\Pi^u x_t-F^s(\Pi^s x_t)\|_X 
   \leq C \|\Pi^u x_{t_1}-F^s(\Pi^s x_{t_1})\|_Xe^{\beta_1 (t-t_1)}, 
          \quad t\in  J
$$
holds true.     
\end{enumerate}
\end{Thm}

The properties (i) through (iii) of Theorem \ref{Theorem D6} can be 
proved in a similar manner to \cite[Theorem 5]{murnag}.  Indeed, let 
$\beta$ be a positive number less than $\alpha$, and $Y_\beta^+$ 
the Banach space $BC^{\beta}({\bb R}^{+};X)$, that is,
$$
  Y_\beta^+
  := BC^\beta({\bb R}^+; X)
  = \big\{y\in C({\bb R}^+;U): 
    \sup_{t\in {\bb R}^+}\|y(t)\|_{X}\,e^{\beta t}<\infty \big\}
$$
with norm $\|y\|_{Y_\beta^+}:=\sup_{t\in {\bb R}^+} \|y(t)\|_X e^{\beta t}$ 
for $y\in Y_\beta^+$. For sufficiently small $r_0>0$ and $\delta>0$ 
one can define a map 
${\cal F}: B_{E^{s}}(r_0)\times B_{Y_\beta^+}(\delta )\to Y_\beta^+$ by 
\begin{align*}
  {\cal F}(\psi,y)(t)
  &:= y(t)-T^s(t)\psi
      -\lim_{n\to \infty}\int_{0}^{t}T^s(t-s)\Pi^{s} \Gamma^n f(y(s))ds  \\
  &\quad  + \lim_{n\to \infty}\int_{t}^{\infty}T^u(t-s)\Pi^{u} \Gamma^n f(y(s))ds 
\end{align*}
for $(\psi,y)\in B_{E^{s}}(r_0)\times B_{Y_\beta^+}(\delta)$ and 
$t\in {\bb R}^+$. In view of the $C^k$-smoothness of the induced map 
$\widetilde f: B_{Y_\beta^+}(\delta ) \to BC({\bb R}^+;{\bb C}^m)$ 
(cf. \cite[Appendix IV]{diek} and \cite{murnag}), the map ${\cal F}$ 
turns out, in contrast to ${\cal F}_{\delta}$, to be of class $C^k$. 
The implicit function theorem (e.g. see \cite [Theorem 5.9 in Chapter I]{lan}) 
then yields the existence of the $C^k$-map 
$\Lambda^s:B_{E^s}(r)\to B_{Y_\beta^+}(\delta)$ satisfying  
${\cal F}(\psi,\Lambda^s(\psi))=0$ for $\psi\in B_{E^s}(r)$, 
$r$ being some positive number with $r\leq r_0$. $\Lambda ^s(\cdot)$ 
plays a similar role to the one $\Lambda _{\ast,\delta}(\cdot)$ does 
in the construction of center manifolds and, thus, 
$F^s:=\Pi^u\circ \mathrm{ev}_0\circ \Lambda ^s$ is the desired 
one satisfying Properties (i), (ii) and (iii); we omit the details.  

Property (iv), which is a result parallel to the result 
\cite[Theorem 13.5.1]{codlev} for ordinary differential equations, can 
be established by similar arguments to Propositions \ref{Proposition D1} 
through \ref{Proposition D6}; so, we will give only a sketch of the proof. 
Given $\psi\in E^s$, the equation for $y\in Y^+_{\beta}$ 
\begin{align*}
   y(t)
   &= T^s(t)\psi+\lim_{n\to \infty}\int_{0}^{t}
      T^s(t-s)\Pi^{s} \Gamma^n f_{\delta}(y(s))ds\\
   &\quad -\lim_{n\to \infty}\int_{t}^{\infty}
          T^u(t-s)\Pi^{u} \Gamma^n f_{\delta}(y(s))ds, \quad t\in {\bb R}^+,
\end{align*}
possesses a unique solution $\Lambda^{\ast}_{\delta}(\psi)$, where 
$\delta \in (0,\delta_1]$, and 
$f_{\delta}$ is a function defined by 
$f_{\delta}(\phi):=\chi(\|\Pi^u\phi\|_X/\delta)\chi(\|\Pi^s\phi\|_X/\delta)f(\phi),\ \phi\in X$, 
 and $\delta_1$ is a number satisfying 
(\ref{286a5}) (cf. Proposition \ref{Proposition D1}). 
Considering a map $F^{\ast}_{\delta}:E^s\to E^u$ defined by 
$$
   F^{\ast}_{\delta}(\psi)
   = -\lim_{n\to \infty}\int_{0}^{\infty}
     T^u(-s)\Pi^{u} \Gamma^n f_{\delta}(\Lambda^{\ast}_{\delta}(\psi)(s))ds, 
     \quad \psi\in E^s
$$
which is indeed an extension of $F^s$, one can verify that for any 
$\tau\in {\bb R}$ and 
$\hat{\phi}\in W_{\delta}^s:=\{\psi+F^{\ast}_{\delta}(\psi):\psi\in E^s\}$, 
the solution $x(t; \tau,\hat{\phi},f_{\delta})$ of ($E_{\delta}$) 
exists on $[\tau,\infty)$, and it satisfies the relation 
$\Pi^ux_t(\tau,\hat{\phi}, f_{\delta})
=F^{\ast}_{\delta}(\Pi^sx_t(\tau,\hat{\phi},f_{\delta}))$ for any 
$t\geq \tau$ (cf. Proposition \ref{Proposition D3}-(iii)). 
Let $x$ be a solution of Eq.($E_{\delta}$) on an interval 
$J:=[t_0,t_1]$, and let $\tau\in J$. Then, putting 
$\hat{x}_{\tau}:= \Pi^s x_{\tau}+F^{\ast}_{\delta}(\Pi^s x_{\tau})$ 
we get the following inequalities: 
\begin{align*}
   &\|\Pi^sx_t-\Pi^sx_t(\tau,\hat{x}_{\tau},f_{\delta})\|_X
   \leq K\int_{\tau}^te^{{\mu}(\sigma-t)}
        \|\Pi^ux_{\sigma}-\Pi^ux_{\sigma}(\tau,\hat{x}_{\tau},f_{\delta})\|_X
        d\sigma, \quad \tau\leq t\leq t_1;\\
   &\|\Pi^sx_t-\Pi^sx_t(\tau,\hat{x}_{\tau},f_{\delta})\|_X
   \leq K\int_{\tau}^te^{{\mu'}(\sigma-t)}
        \|\Pi^ux_{\sigma}-F^{\ast}_{\delta}(\Pi^sx_{\sigma})\|_Xd\sigma, 
        \quad \tau\leq t\leq t_1;\\
   &\|\Pi^sx_t(t_1,\hat{x}_{t_1},f_{\delta})-\Pi^sx_t(\tau,\hat{x}_{\tau},f_{\delta})\|_X
   \leq K\int_{\tau}^{t_1}e^{{\mu'}(\sigma-t)}
        \|\Pi^ux_{\sigma}-F^{\ast}_{\delta}(\Pi^sx_{\sigma})\|_Xd\sigma, 
        \quad t\geq t_1;
\end{align*}
here $K:=CC_1\zeta_{\ast}(\delta),\ \mu:=\alpha-K$ and 
$\mu'=\mu-KL(\delta)$ (cf. Propositions \ref{Proposition D4} and 
\ref{Proposition D5}). Subsequently, utilizing these results and 
repeating the argument similar to Proposition \ref{Proposition D6} 
we can establish the inequality 
$$
   \|\Pi^ux_t-F^{\ast}_{\delta}(\Pi^sx_t)\|_X
   \leq C\|\Pi^ux_{t_1}-F^{\ast}_{\delta}(\Pi^sx_{t_1})\|_Xe^{\beta_1(t-t_1)}, 
        \quad t\in J,
$$
which implies Property (iv); here $\beta_1$ is a (positive) number 
given by $\beta_1:=\alpha\{1-2K/(2\alpha-K-KL(\delta))\}$.

\medskip

Similarly we can establish the existence of local unstable manifolds 
for Eq.($E$). For $r>0$ and $\delta>0$ consider the set 
$$
   W^{u}_{{\rm loc}}(r,\delta)
   = \big\{\phi \in X: \| \Pi^{u}\phi\|_X < r, 
     \enskip \| x_{t}(0,\phi,f)\|_X <\delta, \ t\in {\bb R}^{-}\big\}.
$$
Then we have the following theorem on the existence of $C^k$-smooth 
local unstable manifolds for Eq.($E$); we omit the proof of the theorem. 

\begin{Thm}\label{Theorem 5.2}
Assume that the zero solution of ($E$) is hyperbolic and that 
$f\in C^{k}(X;{\bb C}^m)$ with $f(0)=Df(0)=0$. Then there exist 
positive numbers $r$, $\delta$, and a $C^k$-map $F^u:B_{E^{u}}(r)\to E^{s}$ 
with $F^u(0)=0$, together with an open neighborhood $\Omega_0$ of 
$0$ in $X$, such that the following properties hold: 
\begin{enumerate}
   \item $W^{u}_{{\rm loc}}(r,\delta)={\rm graph}\,F^u$; moreover,  $W^{u}_{{\rm loc}}(r,\delta)$ is 
         tangent to $E^{u}$ at zero.  
   \item To any $\beta\in (0,\alpha)$ there corresponds a positive 
         constant $M$ such that 
$$
   \| x_{t}(0,\phi,f)\|_X
   \leq Me^{\beta t}\|\phi\|_X, 
        \quad t\in {\bb R}^{-}, \enskip \phi\in W^{u}_{{\rm loc}}(r,\delta),
$$
   \item $W^{u}_{{\rm loc}}(r,\delta)$ is locally negatively invariant 
         for ($E$), that is, if $\phi\in W^{u}_{{\rm loc}}(r,\delta)$, 
         we have $x_{\tau}(0,\phi,f)\in W^{u}_{{\rm loc}}(r,\delta)$ 
         for $\tau\in {\bb R}^{-}$ whenever 
         $\Pi^{u}x_{\tau}(0,\phi,f)\in B_{E^{u}}(r)$.    
   \item There exists a positive constant $\beta_1$ with the property 
         that if $x(t)$ is a solution of Eq.($E$) on an interval 
         $J=[t_0,t_1]$ satisfying $x_t\in \Omega_0$ on $J$, then the 
         inequality 
$$
   \|\Pi^s x_t-F^u(\Pi^u x_t)\|_X 
   \leq C \|\Pi^s x_{t_0}-F^u(\Pi^u x_{t_0})\|_Xe^{-\beta_1 (t-t_0)}, 
        \quad t\in  J
$$
holds true.  
\end{enumerate}
\end{Thm}

\begin{Rem} {\rm 
We can also establish the existence and $C^k$-smoothness of 
center-stable/center-unstable manifolds for Eq.($E$) provided that 
$f$ is of class $C^k$. We will, however, omit the statements and 
the proofs of the theorems.  }     
\end{Rem}

\vspace{5mm}
\bibliographystyle{amsplain}

\end{document}